\documentclass[reqno]{amsart}
\newtheorem{thm}{Theorem}[section]
\newtheorem{df}[thm]{Definition}
\newtheorem{prop}[thm]{Proposition}
\newtheorem{lem}[thm]{Lemma}
\newtheorem{rem}[thm]{Remark}
\newtheorem{cor}[thm]{Corollary}
\newtheorem{exm}[thm]{Example}
\usepackage{amsmath,amssymb,amscd,epsfig}
\usepackage{color, graphics}
\numberwithin{equation}{section}

\newcommand{\ld}{\lambda}

\begin{document}
\title[$q$-deformation of Witt-Burnside rings]
{$q$-deformation of Witt-Burnside rings}
\author{YOUNG-TAK OH}
\address{Department of Mathematics \\ Sogang University\\Seoul 121-742, Korea}
\email{ytoh@sogang.ac.kr} \maketitle \baselineskip=12pt


\begin{abstract}
In this paper, we construct a $q$-deformation of the Witt-Burnside
ring of a profinite group over a commutative ring, where $q$ ranges
over the set of integers. When $q=1$, it coincides with the
Witt-Burnside ring introduced by A. Dress and C. Siebeneicher (Adv.
Math. {70} (1988), 87-132). To achieve our goal we first show that
there exists a $q$-deformation of the necklace ring of a profinite
group over a commutative ring. As in the classical case, i.e., the
case $q=1$, $q$-deformed Witt-Burnside rings and necklace rings
always come equipped with inductions and restrictions. We also study
their properties. As a byproduct, we prove a conjecture due to Lenart
(J. Algebra. 199 (1998), 703-732). Finally, we classify
$\mathbb W_G^q$ up to strict natural isomorphism in case where $G$
is an abelian profinite group.
\end{abstract}

\section{Introduction}
\renewcommand{\thefootnote}{}
\footnotetext[1]{%
\renewcommand{\baselinestretch}{1.2}\selectfont
The author gratefully acknowledges support from the following grants: KOSEF Grant \# R01-2003-000-10012-0;
KRF Grant \# 2006-331-C00011.
\hfill \break MSC : 11F03,11F22,17B70.
\hfill \break Keywords: Witt-vectors, Necklace ring, Witt-Burnside
ring, Burnside-Grothendieck ring}

The universal ring of Witt vectors was introduced by Witt and Lang around 1965.
In \cite {La}, they discovered that
there exists a unique covariant functor $\mathbb W$ on the category of commutative rings with identity
characterized by the
following properties:
\begin{enumerate}
\item
As a set, it is $A^{\mathbb N}.$
\item
For any ring homomorphism  $f:A\to B$, the map
$\mathbb W(f):{\bf a}\mapsto (f(a_n))_{n\ge 1}$
is a ring homomorphism
for ${\bf a}=(a_n)_{n\ge 1}$.
\item
The maps $w_m: \mathbb W(A)\to A$ defined by
$${\bf a} \mapsto \sum_{d|m}da_d^{\frac md}
\text{  for }{\bf a}=(a_n)_{n\ge 1}$$
are ring homomorphisms (\cite{R,H,La})
\end{enumerate}
Here, $\mathbb W(A)$ is called the {\it universal ring of Witt vectors over $A$}.

The theory of Witt vectors has been developed with deep connection with necklace rings,
which are algebraic objects arising from the combinatorics of necklaces.
To be more precise, in \cite {MR}, Metropolis and Rota introduced a covariant functor $Nr$
on the category of commutative rings with identity.
In this case, $Nr(A)$ is called the {\it necklace ring over $A$}.
They then showed that
\begin{equation*}\label{isom1}
\mathbb W(\mathbb Z)\cong Nr(\mathbb Z).
\end{equation*}
To do this they studied the combinatorics of necklaces extensively and described Witt vectors
in terms of {\it strings}, which are infinite matrices satisfying suitable conditions.

On the other hand, in \cite{DS}, Dress and Siebeneicher observed that ${Nr}(\mathbb Z)$
is naturally isomorphic to the Burnside-Grothendieck ring
$\hat \Omega(C)$ of isomorphism classes
of almost finite cyclic sets.
Here, $C$ denotes the multiplicative infinite cyclic group.
Motivated by this fact,
they found that given any profinite group, say $G$,
there exists a unique covariant functor $\mathbb W_G$
on the category of commutative rings with identity,
satisfying
\begin{equation*}\label{isom2}
{\rm i)}\,\,\, \mathbb W_G(\mathbb Z)\cong \hat \Omega(G)\quad \text{ and } \quad
{\rm ii)}\,\,\,\mathbb W_{\hat C}=\mathbb W.
\end{equation*}
Here, $\hat \Omega(G)$ denotes the Burnside-Grothendieck ring of $G$
of isomorphism classes of almost finite $G$-spaces, and $\hat C$ the profinite completion of $C$.
In this case, $\mathbb W_G(A)$ is called the {\it Witt-Burnside ring of $G$ over $A$}.

Recently, in \cite{O,O2}, it was shown that given any profinite group $G$,
there exists a unique covariant functor $Nr_G$ from the category of special $\ld$-rings
to the category of commutative rings with identity, satisfying
$$\mathbb W_G(A)\cong Nr_G(A),$$
where $A$ is any special $\ld$-ring.
Obviously, $Nr_{\hat C}$ is not equivalent to $Nr$ since domains are different.
This led us to introduce $\widehat {Nr}_G$,
which is a unique covariant functor on the category of commutative rings with identity, satisfying
$${\rm i)}\,\,\, \widehat {Nr}_G(\mathbb Z)\cong \hat \Omega(G)
\quad \text{ and } \quad {\rm ii)}\,\,\,\widehat {Nr}_{\hat C}=Nr.$$
In view of ii), $\widehat {Nr}_{G}$ may be regarded as a group-theoretic generalization of $Nr$.
Furthermore, viewed as a functor
from the category of binomial rings to the category of commutative rings with identity,
$Nr_G$ and $\widehat {Nr}_{G}$ turn out to be same (\cite{O3,O2}).

The theory of Witt vectors arises in the context of formal group laws, too (\cite{H}).
It is well-known that for every one-dimensional formal group law
$F(X,Y)$ over a torsion-free ring $B$, there exists a unique covariant
functor $\mathbb W^F$ from the category
of $B$-algebras to the category of abelian groups.
Related to this functor, it is quite surprising that if $F(X,Y)$
is given by $X+Y-qXY,\,q\in \mathbb Z$, then the value of $\mathbb
W^F$, usually denoted by $\mathbb W^q$, at each object has the structure of a commutative ring.
This phenomenon was first observed by Lenart (\cite{L})
for commutative torsion-free rings with identity.
Later, in \cite {O3}, it was shown that his observation makes sense for arbitrary commutative rings.
Consequently, one can derive a covariant functor $\mathbb W^q$ from the
category of commutative rings to itself in this case.
In particular, if $q=1$, then it coincides with $\mathbb W$.
In this sense, we may regard $\mathbb W^q$ as a $q$-deformation of $\mathbb W$.

Since $\mathbb W=\mathbb W_{\hat C}$, it would be quite natural to expect
that the same phenomenon also
happens to the functor $\mathbb W_G$ for arbitrary profinite groups.
Recall that the key ingredient in proving the existence of $\mathbb W_G$
is that the polynomials determining the ring structure of $\mathbb W_G(\mathbb Z)$ have integer coefficients.
These polynomials were utilized to endow $\mathbb W_G(A)$ with the ring structure
for arbitrary commutative rings with identity.
In contrast, the polynomials arising from the construction of
$\mathbb W_G^q$, the $q$-deformation of $\mathbb W_G$,
are in $\mathbb Q[q]$.
The most significant and difficult step in our construction will be to show
that these polynomials
take integer values for integer arguments.
To do this we first introduce $Nr_G^q$, the $q$-deformation of $Nr_G$,
and then define $q$-inductions and $q$-restrictions.
Exploiting $q$-inductions and $q$-restrictions,
we construct an isomorphism, called $q$-Teichm\"{u}ller map.
With this preparation, we finally show that
Dress and Siebeneicher's proof can be applied to our case.

This paper is organized as follows.
In Section \ref{preliminaries},
we introduce the functors $\mathbb W_G$, $Nr_G$ and $\widehat {Nr}_G$.
In Section \ref{orbit},
we introduce orbit-sum polynomials associated with a profinite group, which play the same role as necklace polynomials
in the theory of Witt vectors.
In Section \ref{Necklace},
$Nr_G^q$ and $\widehat {Nr}_G^q$ will be defined as $q$ ranges over the set of integers,
and in Section \ref{on indiction and restriction},
natural transformations called $q$-induction and $q$-restriction will be defined.
Main results will appear in Section \ref{Witt-Burn}.
$\mathbb W_G^q$ will be defined as $q$ ranges over the set of integers, and
many structural properties will be investigated.
Section \ref{Proofs} is devoted to proving lemmas and theorems which are stated without proof.
In particular, we prove a conjecture due to Lenart(\cite {L})
concerning $q$-restriction and $q$-necklace polynomials.
In the final section, we show that $\mathbb W_G^q$ is classified up to strict natural isomorphism
by ${\rm D}^{\rm pr}(q) \cap {\rm D}^{\rm pr}(G)$ as $q$ ranges over the set of integers.
Here,
${\rm D}^{\rm pr}(q)$ denotes the set of prime divisors of $q$, and
${\rm D}^{\rm pr}(G)$ the set of prime divisors of each of
$$\{(G:U)\,\,|\,\,U \text{ is an open subgroup of } G\}.$$


\section{The Witt- Burnside ring and the necklace ring of a profinite group}
\label{preliminaries}

In this section, we review prerequisites on
the Witt- Burnside ring and the necklace ring
of a profinite group. For more information see \cite{O,O2}.
Unless otherwise stated, the rings we consider
will be commutative, but not necessarily unital.

Let $G$ be an arbitrary profinite group.
For any $G$-space $X$ and any subgroup $U$
of $G$ define $\varphi_U(X)$ to be the cardinality
of the set $X^U$ of $U$-invariant
elements of $X$ and let $G/U$ denote the $G$-space of
left cosets of $U$ in $G.$
For two subgroups $U,\,V$ of $G$, we say that
$U$ is subconjugate to $V$ if $U$ is a subgroup
of some conjugates of $V$.
This is a partial order on the set of
the conjugacy classes of open subgroups
of $G$, and will be denoted by $[V]\preceq [U]$.
Fix an enumeration of this poset
satisfying the condition:\\
\centerline
{\it If $[V]\preceq [U]$, then $[V] \text{ precedes }[U].$}\\
By abuse of notation we denote this poset by $\mathcal O(G)$.
For example, if $G$ is abelian${}^*$,
then $\mathcal O(G)$ is just the set of open subgroups of $G$ subject to
$$V \preceq U \Longleftrightarrow
U \subseteq  V .$$
\footnote{${}^*$ In case where $G$ is abelian, we omit the bracket notation.}

Given a commutative ring $A$ and a profinite group $G$, we define the
{\it ghost ring of $G$ over $A$}, denoted by ${\rm Gh}(G,A)$,
to be the commutative ring
$$A^{\mathcal O(G)}\,\left(:=\underset{i\in \mathcal O(G)}{{\prod}} A\right)$$
whose addition and multiplication are defined componentwise.
\begin{df}{\rm (\cite{DS2})}
Given a profinite group $G$, $\mathbb W_G$ is a unique covariant functor
from the category of commutative rings with identity to itself
characterized as follows.
\begin{enumerate}
\item
As a set
$$\mathbb W_G(A)=A^{\mathcal O(G)}.$$
\item
For every ring homomorphism $f:A\to B$
and every ${\bf x} \in  \mathbb W_G(A)$ one has
$$\mathbb W_G(f)({\bf x})=f\circ {\bf x}.$$
\item
The map,
\begin{equation*}
\Phi:  \mathbb W_G(A) \to {\rm Gh}(G,A),\quad
{\bf x} \mapsto \left(\sum_{[G]\preceq [V] \preceq [U]}
\varphi_U(G/V)\cdot {\bf x}([V])^{(V:U)}\right)_{[U]\in \mathcal O(G)},
\end{equation*}
is a ring homomorphism.
Here, $(V:U)$ represents $(G:U)/(G:V)$.
\end{enumerate}
In this case, $\mathbb W_G(A)$ is called
the {\it Witt-Burnside ring of $G$ over $A$}.
\end{df}
\begin{exm}\hfill
{\rm

(a) If $G=\hat C$,
the profinite completion of the multiplicative infinite cyclic group $C$,
then the conjugacy classes of open subgroups are parametrized
naturally by their index in $\hat C.$
Thus, $\mathbb W_G$ coincides with $\mathbb W$, the functor of the ring of Witt vectors
due to Witt and Lang (\cite{Di,H}).

(b) If $G=\hat C_p$, the pro-$p$-completion of the infinite cyclic group, then
$\mathbb W_G$ coincides with $\mathbb W_p$, the functor of the ring of $p$-typical Witt vectors
due to Witt (\cite{H, WI}).

(c) For a positive integer $n$ we denote by
${\bf n}$ the set of divisors of $n$.
If $G$ is the finite cyclic group of order $n$, then $\mathbb W_G$ coincides with
the functor of the ring of $\bf n$-nested Witt vectors (\cite{R}).
}\end{exm}

Let us introduce a functor which is naturally equivalent to $\mathbb W_G$ on a category of special $\ld$-rings.
References for special $\ld$-rings are \cite{AD,DK,O,O3,O2}.
Define a $\mathcal O(G)\times \mathcal O(G)$ matrix
$\tilde \zeta_{G}$ by
\begin{equation*}
\tilde \zeta_G([V],[W])=
\begin{cases}
\varphi_V(G/W)\Psi^{(W:V)} &\text{ if }[W]\preceq [V],\\
0 & \text{ otherwise.}
\end{cases}
\end{equation*}
Here, the notation $\Psi^{(W:V)}$ represents the $(W\!\!:\!\!V)$-th Adams operation.
\begin{df}{\rm (\cite{O3,O2})}
Given a profinite group $G$, ${Nr}_G$ is a unique covariant functor
from the category of special $\ld$-rings to the category of commutative rings with identity
characterized as follows.
\begin{enumerate}
\item
As a set
$$Nr_G(A)=A^{\mathcal O(G)}.$$
\item
For every special $\ld$-ring homomorphism $f:A\to B$ and
every $\alpha \in  Nr_G(A)$
one has
$$Nr_G(f)(\alpha)=f\circ \alpha.$$
\item
The map,
\begin{equation*}
\tilde\varphi: Nr_G(A) \to {\rm Gh}(G,A),\quad
{\bf x} \mapsto \tilde \zeta_G \,{\bf x},
\end{equation*}
is a ring homomorphism. Here, we understand ${\bf x}$
as a $1\times \mathcal O(G)$ column vector.
\end{enumerate}
\end{df}
While the addition in $Nr_G(A)$ is defined componentwise,
the multiplication in $Nr_G(A)$ is somewhat complicated.
Given ${\bf x}, {\bf y} \in Nr_G(A)$ and $[U]\in \mathcal O(G)$,
the $[U]$-th component of ${\bf x\cdot y}$ is given by
\begin{equation}\label{mu1}
\underset{[V],[W] \in\mathcal O(G)}
{{\sum}}\sum_{VgW\subseteq G \atop [Z(g,V,W)]=[U]}
\Psi^{(V:Z(g,V,W))}({\bf x}([V]))\,\Psi^{(W:Z(g,V,W))}({\bf y}([W])),
\end{equation}
where $Z(g,V,W)$ denotes $V\cap gWg^{-1}.$

\begin{exm}
{\rm
If $G$ is abelian, Eq. \eqref{mu1} is reduced to the following simple form
\begin{equation*}
\underset{V,W \in\mathcal O(G)}
{{\sum}}\sum_{V\cap W=U}(G:V+W)
\Psi^{(V:U))}({\bf x}(V))\,\Psi^{(W:U)}({\bf y}(W)).
\end{equation*}
In particular, if $G=\hat C$, then
\begin{equation*}
{\bf x}\cdot {\bf y}(n)=\sum_{[i,j]=n}(i,j)
\Psi^{(\frac nd)}({\bf x}(d))\,\Psi^{\frac nd}({\bf y}(d)).
\end{equation*}
Here, $[i,j]$ represents the least common multiple
and $(i,j)$ the greatest common divisor of $i$ and $j.$
}\end{exm}

\begin{thm}\label{well-def of neck functors} {\rm (\cite{O2})}
Viewed as a functor from the category of special $\ld$-rings
to the category of commutative rings with identity,
$\mathbb W_G$ is naturally isomorphic to $Nr_G$.
\end{thm}

\begin{rem}{\rm
A commutative ring $A$ with identity may have two different special $\ld$-ring structures.
Each of them will produce two different commutative rings, say $Nr_{G,1}(A)$ and $Nr_{G,2}(A).$
However, Theorem \ref {well-def of neck functors} says that they are isomorphic.
}\end{rem}

Define a $\mathcal O(G)\times \mathcal O(G)$ matrix
$\zeta_{G}$${}^\dag$
\footnote{${}^\dag$ In the literature this matrix is called the {\it Burnside matrix of $G$}.}
by
\begin{equation*}
\zeta_{G}([V],[W])=\varphi_V(G/W)
\end{equation*}
for all $[V],[W]\in \mathcal O(G).$ With this notation, we have
\begin{df}{\rm (\cite{O3})}
Given a profinite group $G$, ${\widehat {Nr}}_G$ is a unique covariant functor
from the category of commutative rings with identity to itself
characterized as follows.
\begin{enumerate}
\item
As a set
$${\widehat {Nr}}_G(A)=A^{\mathcal O(G)}.$$
\item
For every ring homomorphism $f:A\to B$ and
every $\alpha \in  {\widehat {Nr}}_G(A)$
one has
$${\widehat {Nr}}_G(f)(\alpha)=f\circ \alpha.$$
\item
The map,
\begin{equation*}\hat \varphi: {\widehat {Nr}}_G(A) \to {\rm Gh}(G,A),\quad
{\bf x} \mapsto \zeta_G \,{\bf x},
\end{equation*}
is a ring homomorphism.
Here, we understand ${\bf x}$
as a $1\times \mathcal O(G)$ column vector.
\end{enumerate}
\end{df}
Note that ${\widehat {Nr}}_G(\mathbb Z)$ is isomorphic to the Burnside-Grothendieck ring of
almost finite $G$-spaces.
The addition in ${\widehat {Nr}}_G(A)$ is defined componentwise.
On the other hand, the multiplication is defined as follows:
Given ${\bf x}, {\bf y} \in {\widehat {Nr}}_G(A)$ and $[U]\in \mathcal O(G)$,
\begin{equation}\label{mu2}
{\bf x\cdot y}([U])=\underset{[V],[W] \in\mathcal O(G)}
{{\sum}}\sum_{VgW\subseteq G \atop [Z(g,V,W)]=[U]}
{\bf x}([V])\,{\bf y}([W]).
\end{equation}
\begin{exm}
{\rm
If $G$ is abelian, Eq. \eqref{mu2} is reduced to
\begin{equation*}
\underset{V,W \in\mathcal O(G)}
{{\sum}}\sum_{V\cap W=U}(G:V+W)
{\bf x}(V)\,{\bf y}(W).
\end{equation*}
In particular, if $G=\hat C$, then
\begin{equation*}
{\bf x}\cdot {\bf y}(n)=\sum_{[i,j]=n}(i,j){\bf x}(d)\,{\bf y}(d),\quad n\ge 1.
\end{equation*}
Indeed, ${\widehat {Nr}}_{\hat C}$ was first introduced by Metropolis and Rota in order to describe Witt vectors (\cite{MR}).
}\end{exm}

Recall that a special $\ld$-ring in which $\Psi^n=id$ for all $n\ge 1$
is called a {\it binomial ring}.
The reason why binomial rings are important is that
if $A$ is a binomial ring, then $Nr_G(A)$ coincides with $\widehat{Nr}_G(A)$.
This follows from Eq. \eqref{mu1} and Eq. \eqref{mu2} immediately.
Thus, viewed as a functor from the category of binomial rings to the category of commutative rings with identity,
$Nr_G$ and ${\widehat {Nr}}_G$ will denote the same functor.
Consequently we can conclude that
\begin{enumerate}
\item
$Nr_G\cong \mathbb W_G$ as a functor on the category of special $\ld$-rings, and
\item
$Nr_G=\widehat {Nr}_G$ as a functor on the category of binomial rings.
\end{enumerate}
\begin{rem}${}^*$ \label{footnote}
\footnote{${}^*$ The author is indebted to the referee for Remark \ref{footnote}.}
{\rm
Since we have necklace ring functors for both commutative rings with identity and special $\ld$-rings,
it would be quite natural to expect Witt-Burnside ring functors
for both commutative rings with identity and special $\ld$-rings, which has not, to the author's knowledge, been known yet.
So, it seems quite challenging to construct a functor, say ${\widehat {\mathbb W}}_G$, satisfying
\begin{enumerate}
\item
${\widehat {\mathbb W}}_G\cong \widehat {Nr}_G$ on the category of special $\ld$-rings, and
\item
${\widehat {\mathbb W}}_G=\mathbb W_G$ on the category of binomial rings.
\end{enumerate}
}\end{rem}
\section{Orbit-sum polynomials associated with profinite groups}
\label{orbit}

Let $G$ be a profinite group and
$X$ be an alphabet, that is,
a set of commuting variables $\{x_1,x_2,\cdots,x_m\}$.
Let us consider a discrete $G$-space (with respect to the discrete topology)
$$\mathbb Z/q\mathbb Z \times X=\{(\bar c,x): \bar c \in \mathbb Z/q\mathbb Z, x \in X \}$$
with trivial $G$-action.
Denote by $\mathcal C(G,\mathbb Z/q\mathbb Z \times X)$
the set of all continuous maps from $G$ to
$(\mathbb Z/q\mathbb Z \times X)$.
It is well known that this set becomes
a $G$-space with regard to the compact-open topology
via the following standard $G$-action
$$(g\cdot f)(a)=f(g^{-1}\cdot a)\quad
\text{ for all } a,g\in G$$
(refer to \cite{DS2}).
For every open subgroup $V$ of $G$,
the set of $V$-invariant elements in
$\mathcal C(G,\mathbb Z/q\mathbb Z \times X)$
coincides with
$${\rm Hom}_V(G,\mathbb Z/q\mathbb Z \times X),$$
the set of continuous $V$-maps from $G$ to $(\mathbb Z/q\mathbb Z \times X)$.
If $f\in {\rm Hom}_V(G,\mathbb Z/q\mathbb Z \times X)$, then
it is contained in an orbit isomorphic to $G/W$ for some $W$
to which $V$ is subconjugate.
On the other hand, $(\mathbb Z/q\mathbb Z \times X)$
becomes a $\mathbb Z/q\mathbb Z$-set
via the action on the first component.
Then $\mathcal C(G,\mathbb Z/q\mathbb Z \times X)$
also becomes a $\mathbb Z/q\mathbb Z$-set
via the action
\begin{equation*}
(\bar c\cdot f)(x)= (\bar c \cdot (\pi_1\circ f)(x), (\pi_2\circ f) (x)), \quad (\bar c\in \mathbb Z/q\mathbb Z, \,\,x\in G).
\end{equation*}
Here, $\pi_i$, $(i=1,2),$ represents the projection to the $i$-th component.
With regard to $\mathbb Z/q\mathbb Z$-action we denote by
$$\mathcal C(G,\mathbb Z/q\mathbb Z \times X)/\sim$$
the set of the equivalence classes of
all continuous maps from $G$ to $(\mathbb Z/q\mathbb Z \times X)$.
In the natural way it becomes a $G$-set.
Decompose $\mathcal C(G,\mathbb Z/q\mathbb Z \times X)/\sim \,$
into disjoint $G$-orbits and then consider its disjoint union,
say,
\begin{equation}\label{disjoint union**}
\underset{h}{\dot\bigcup}
\,Gh.
\end{equation}
Here, $h$ runs through a system of representatives
of this decomposition.
Let $G_h$ be the isotropy subgroup of $h$.
Write
$$G=\underset{1\le i \le (G:G_h)}{\dot\bigcup}\,\,G_hw_i,$$
where $\{w_i:1\le i \le (G:G_h)\}$ is
a fixed set of right coset representatives.
With the above notation, we define $[h]$ by the polynomial in $x_1,x_2,\cdots,x_m$ given by
$$\prod_{i=1}^{(G:G_h)}(\pi_2\circ h)(w_i).$$
It is not difficult to show that it is well defined, that is, it does not depend on the choice of coset representatives.
\begin{df}
Let $h\in {\rm Hom}_V(G,\mathbb Z/q\mathbb Z \times X)/\sim$ and $W$ be an open subgroup of $G$.
We say that $h$ has a period $W$ if
\begin{enumerate}
\item
$G_h \subseteq W,$
\item
There exists an element $h_W \in {\rm Hom}_W(G,\mathbb Z/q\mathbb Z \times X)/\sim$
and $s_j \in \mathbb Z/q\mathbb Z,\,1\le j \le (W:G_h),$ such that
$$h(t_jw_i)=s_j\cdot h_W(w_i)$$
for $1\le i \le (G:W),\,\,1\le j \le (W:G_h)$.
Here, $\{w_i : 1\le i\le (G:W)\}$ is a set of
right-coset representatives of $W$ in $G$
and $\{t_j : 1\le j\le (W:G_h)\}$
a set of right-coset representatives of $G_h$ in $W$.
\end{enumerate}
If $h$ has a period $G_h$,
then it is called {\it aperiodic}.
\end{df}

\begin{df}\label{def o fq-necklace polyno}
Let $q$ be any positive integer and $G$ be a profinite group.

{\rm (a)}
Given an open subgroup $V$ of $G$, we define $M^q_G(X,V)$, called the {\it orbit-sum} polynomial of $V$,
by
the polynomial in $x_1,x_2,\cdots, x_m$ given by
$$\sum_h\,[h].$$
Here, $h$ ranges over the set of aperiodic representatives
in the decomposition \eqref{disjoint union**} such that
$Gh$ is isomorphic to $G/V$.

{\rm (b)}
Let $f\in {\rm Hom}_V(G,\mathbb Z/q\mathbb Z \times X)/\sim$.
Then the {\it weight of $f$ over $V$}, denoted by ${\rm wt}_V(f)$,
is defined by
\begin{equation*}
\prod_{i=1}^{(G:V)}(\pi_2\circ f)(v_i).
\end{equation*}
Here, $v_i$'s range over a set of
right-coset representatives of $V$ in $G$.
\end{df}

Define $\Psi^n$, the $n$-th Adams operation, on
$\mathbb Q\,[\,x_i\,:\,1\le i\le m\,]$ as follows:
\begin{align*}
&\Psi^n(x_i)=x_i^n,\quad 1\le i\le m, \,n\ge 1\\
&\Psi^n(c)=c, \quad c\in \mathbb Q.
\end{align*}
The following lemma illustrates an intrinsic relation between orbit-sum polynomials and
weights of functions.
\begin{lem}\label{relation of hom and necklace*}
Let $q$ be any positive integer and $V$ be an open subgroup of $G$.
Then the following identity holds.
\begin{equation}\label{first main formula}
\sum_{f\in {\rm Hom}_V(G,\mathbb Z/q\mathbb Z \times X)/\sim}
{\rm wt}_V(f)
=\sum_{[W]\preceq [V]}\varphi_V(G/W)q^{(W:V)-1} \Psi^{(W:V)}(M^q_G(X,W)).
\end{equation}
\end{lem}

\noindent{\bf Proof.}
Decompose $\mathcal C(G,\mathbb Z/q\mathbb Z \times X)/\sim$
into disjoint $G$-orbits
and consider its union, say,
$$\underset{h}{\dot\bigcup}\,G h.$$
Observe that $V$-invariant maps exist only in the orbits
such that $Gh$ are isomorphic to $G/W$ for some $W$
to which $V$ is subconjugate.
For each aperiodic function $h$ such that $Gh$ is isomorphic to $G/W$,
note that there exist $q^{(W:V)-1}$-number of maps with period $h$
in ${\rm Hom}_V(G,\mathbb Z/q\mathbb Z \times X)/\sim$.
On the other hand, every element in
${\rm Hom}_V(G,\mathbb Z/q\mathbb Z \times X)/\sim$
arises in this way.
Since the number of all $V$-invariant elements in the orbit $Gh$ is
given by $\varphi_V(G/W)$ by definition and
$${\rm wt}_V(f)=[f]^{(W:V)}$$
for all $V$-invariant functions $f \in Gh$,
the desired result follows.
\qed
\begin{exm}{\rm
Note that
\begin{equation}\label{exam of neck and word}
\sum_{f\in {\rm Hom}_V(G,\mathbb Z/q\mathbb Z \times X)/\sim}
{\rm wt}_V(f)=q^{(G:V)-1}(x_1+\cdots+x_m)^{(G:V)}.
\end{equation}
Thus, if $G=\hat C$ and $V=\hat C^n$ (a unique open subgroup of $\hat C$ of index $n$),
then Eq. \eqref{first main formula} is reduced to
\begin{equation*}
q^{n-1}(x_1+\cdots+x_m)^n
=\sum_{d|n}dq^{\frac nd-1} \Psi^{\frac nd}(M^q_G(X,\hat C^d)).
\end{equation*}
}\end{exm}
From now on, we let $q$ be an indeterminate.
Let us define
a $\mathcal O(G)\times \mathcal O(G)$ matrix
$\tilde \zeta^q_{G}$ by
\begin{equation}\label{def of q-matrix }
\tilde \zeta^q_G([V],[W])=
\begin{cases}
\varphi_V(G/W)\,q^{(W:V)-1}\,\Psi^{(W:V)}
& \text{ if }[W]\preceq [V]\\
0 & \text { otherwise.}
\end{cases}
\end{equation}
We also define a $\mathcal O(G)\times \mathcal O(G)$ matrix $\zeta^q_G$ by
\begin{equation}\label{def of q-matrix: standard}
\zeta_{G}^q([V],[W])=
\begin{cases}
\varphi_V(G/W)\,q^{(W:V)-1}& \text{ if }[W]\preceq [V]\\
0 & \text { otherwise.}
\end{cases}
\end{equation}
From the fact that $\varphi_V(G/W)=0$ unless $[W] \preceq [V]$
it follows that $\tilde \zeta_G^q$ is a lower-triangular matrix
with the diagonal elements $(N_G(V):V)\cdot {\rm Id}$,
the index of $V$ in its normalizer $N_G(V)$ in $G$, for $[V]\in \mathcal O(G)$.
Therefore,
$\tilde \zeta_G^q$ is invertible over a ring $A$ with identity if and only if
$(N_G(V):V)\cdot 1$ is a unit in $A$ for all $[V]\in \mathcal O(G)$.
Assume that the base ring is $\mathbb Q[q]$.
Let $\tilde \mu^q_G$ (resp. $\mu^q_G$) be the inverse of $\tilde\zeta_G^q$
(resp. $\zeta_{G}^q$).
Let us investigate how to compute the inverse of $\tilde \mu^q_G$ and $\mu^q_G$.
First,
we let $\mathcal O(G,U)$ be the subset of $\mathcal O(G)$
consisting of the elements satisfying the condition $[V]\preceq [U]$.
Set
$$\tilde\mu_{G,U}^q:={(\tilde \zeta_{G,U}^q)^{-1}} \quad (\text{resp, }\mu_{G,U}^q:={(\zeta_{G,U}^q)^{-1}}),$$
where $\tilde \zeta_{G,U}^q$ (resp. $\zeta_{G,U}^q$) is the matrix obtained from $\tilde \zeta_{G}^q$ (resp. $\zeta_{G}^q$)
by restricting the index to $\mathcal O(G,U)$.
It is not difficult to show that for all $[V]\in \mathcal O(G)$
$$\tilde\mu_G^q([V],[W])
=\begin{cases}
\tilde \mu_{G,V}^q([V],[W])& \text{ if }[W] \in \mathcal O(G,V),\\
0&\text{ otherwise.}
\end{cases}$$
Similarly, one can show that
$$\mu_G([V],[W])
=\begin{cases}
\mu_{G,V}([V],[W])& \text{ if }[W]\in \mathcal O(G,V),\\
0&\text{ otherwise.}
\end{cases}$$
The next lemma shows the relation between $\tilde \mu_G^q$ and $\mu_G^q$.
\begin{lem}\label{q-simpler form of matrix}
For every $[W],[V] \in \mathcal O(G)$ satisfying $[W] \preceq [V]$, we have
$$\tilde \mu^q_G([V],[W])=\mu^q_G([V],[W])\Psi^{(W:V)}.$$
\end{lem}
\noindent{\bf Proof.}
Consider the following linear system
$${\bf y}=\tilde \zeta_{G,V} \,{\bf x}$$
for
${\bf x}, {\bf y}\in \prod_{i\in \mathcal O(G,V)}\mathbb Q[x_i:1\le i\le m].$
This linear system can be rewritten as
$$\tilde {\bf y}=\zeta_{G,V} \,{\bf x},$$
where
$$\tilde {\bf y}([W])=\Psi^{(W:V)}({\bf y}([W]))
\quad \text{ for all }[W]\in \mathcal O(G,V).$$
Thus, we have the following equation:
$${\bf x}
=\tilde \mu_{G,V} \,{\bf y}
=\mu_{G,V} \,\tilde {\bf y}.$$
Computing the $[V]$-th component from each side, we obtain
$$\sum_{[W]\in \mathcal O(G,V)}
\tilde \mu^q_{G,V}([V],[W]){\bf y}([W])
=\sum_{[W]\in \mathcal O(G,V)}\mu^q_{G,V}([V],[W])\Psi^{(W:V)}({\bf y}([W])).$$
This implies
$$\tilde \mu^q_{G,V}([V],[W])=\mu^q_{G,V}([V],[W])\Psi^{(W:V)}.$$
So we are done.
\qed
\begin{thm}\label{expression of q-neck associ profinite group}
Let $X=\{x_1,\cdots,x_m\}$ be an alphabet and $q$ be an indeterminate.
Then, for a profinite group $G$ and an open subgroup $V$ of $G$, we have
\begin{equation}\label{cloased form: q-necklacr poly}
M_G^q(X,V)=\sum_{[W]\preceq [V]}\mu^q_{G}([V],[W])\,
\,q^{(G:W)-1}\,(x_1^{(W:V)}+\cdots+x_m^{(W:V)})^{(G:W)}\,.
\end{equation}
\end{thm}
\noindent{\bf Proof.}
Recall that
$$
\sum_{f\in {\rm Hom}_V(G,\mathbb Z/q\mathbb Z \times X)/\sim}
{\rm wt}_V(f)=q^{(G:V)-1}(x_1+\cdots+x_m)^{(G:V)}$$
(see Eq. \eqref{exam of neck and word}).
From Lemma \ref{relation of hom and necklace*} it follows that
$$\tilde \zeta^q_G
\begin{pmatrix}\vdots\\ M_G^q(X,V)\\ \vdots\end{pmatrix}
=\begin{pmatrix}\vdots\\q^{(G:V)-1}(x_1+\cdots+x_m)^{(G:V)}
\\ \vdots\end{pmatrix}.$$
Taking $\tilde \mu^q_G$ on both sides
and then applying Lemma \ref{q-simpler form of matrix}
one can deduce the desired result.
\qed
\vskip 3mm

Specializing $x_i$ into $1$ for all $1 \le i \le m$,
Eq. \eqref{cloased form: q-necklacr poly} is reduced to
\begin{equation*}
M_G^q(m,V)=\sum_{[W]\preceq [V]}\mu^q_{G}([V],[W])\,q^{(G:W)-1}\,m^{(G:W)}\,
\end{equation*}
which represents the number of aperiodic representatives $h$
such that $Gh$ is isomorphic to $G/V$, in the decomposition \eqref{disjoint union**}.
Replacing $m$ by an indeterminate $x$, we can obtain a polynomial in $x$ and $q$, denoted by $M_G^q(x,V)$.
\begin{lem}\label{criterion of integrality}
Suppose that a polynomial $f(x)\in \mathbb Q[x]$ takes integer values for
all positive integers.
Then it takes integer values for all integers.
\end{lem}

\noindent{\bf Proof.}
Reducing to a common denominator we can write
$f(x)=\frac {g(x)}{c}$
for some $g(x)\in \mathbb Z[x]$ and $c\in \mathbb N$.
Suppose the assertion is false, then there must exist
the largest integer $m_0$ such that
$c$ does not divide $g(m_0)$.
Then $c$ divides $g(m_0+c)$ by the maximal condition of $m_0$.
But, since $g(x+c)=g(x)+c\cdot h(x)$ for some $h(x)\in \mathbb Z[x]$,
$c|g(m_0+c)$ implies that $c|g(m_0)$.
This is a contradiction.
\qed
\vskip 3mm
Lemma \ref{criterion of integrality} says that $M_G^q(x,V)$ is a numerical polynomial in $q$ and $x$,
that is, it takes integer values for all positive integers $q$ and $x$.
\begin{rem}{\rm
The orbit-sum polynomial $M_G^q(X,V)$ is a symmetric polynomial in $x_i$'s.
Hence, it can be written as a polynomial with integral coefficients
in the elementary symmetric functions $e_i(x_1,\cdots,x_m),\,1\le i \le m$.
On the other hand, all the coefficients are contained in $\mathbb Q[q]$.
So, we can conclude that all of them are numerical polynomials in $q$.
}\end{rem}

In case where $G$ is abelian, we often prefer to using
the matrix
${\bf {\tilde\zeta}_{G}^q}$ instead of $\tilde\zeta_G^q$, which is defined by
\begin{equation*}
{\bf \tilde\zeta_{G}^q}([V],[W])=
\begin{cases}
q^{(W:V)-1}\Psi^{(W:V)} &\text{ if } V\subseteq W,\\
0 & \text{ otherwise.}
\end{cases}
\end{equation*}
Let ${\bf \tilde \mu_G^q}$
be the inverse of ${\bf \tilde \zeta_G^q}$.
Then it is easy to show that
$$\tilde \mu_{G}^q(V,W)=\frac {1}{(G:V)}\,\,{\bf \tilde \mu_{G}^q}(V,W).$$
Similarly, if we define
${\bf \zeta_{G}^q}$ by
\begin{equation*}
{\bf \zeta_{G}^q}([V],[W])=
\begin{cases}
q^{(W:V)-1} &\text{ if } V\subseteq W,\\
0 & \text{ otherwise}
\end{cases}
\end{equation*}
and let ${\bf \mu_{G}^q}$ be the inverse of ${\bf \zeta_{G}^q}$,
then it holds
$$\mu_{G}^q(V,W)=\frac {1}{(G:V)}\,\,{\bf \mu_{G}^q}(V,W).$$
Hence, Eq. \eqref{cloased form: q-necklacr poly} is reduced to
\begin{equation}\label{simpler form of generalized necklace}
M_G^q(X,V)=\frac {1}{(G:V)}\sum_{V \subseteq W}{\bf
\mu_{G}^q}(V,W)\,
q^{(G:W)-1}\,(x_1^{(W:V)}+\cdots+x_m^{(W:V)})^{(G:W)}\,.
\end{equation}
\vskip 3mm
In the following, we would like to investigate $M_G^q(x,n)$ in detail in case where
$G=\hat C$, the profinite completion of infinite cyclic group.
Let $\hat C^n$ be a unique open subgroup of $\hat C$ of index $n$.
The polynomial $M_{\hat C}^q(x,\hat C^n),\,(n \in \mathbb N$), first appeared in  \cite {L}
in the context of formal group laws $F_q(X,Y)=X+Y-qXY$, $(q\in \mathbb Z)$.
It was shown that they have a nice combinatorial description
in terms of so called {\it aperiodic $q$-words}.
Let us recall their definition briefly.
First, consider the free monoid $(\mathbb Z/q\mathbb Z \times X)^*$ generated by
$\mathbb Z/q\mathbb Z \times X$.
Then, we define a $\mathbb Z/q\mathbb Z$-action
on $(\mathbb Z/q\mathbb Z \times X)^*$
by letting the generator of $(\mathbb Z/q\mathbb Z)$
act as
\begin{equation*}
\left((r_1,a_1),\cdots,(r_n,a_n)\right)\mapsto
\left((r_1+1,a_1),\cdots,(r_n+1,a_n)\right).
\end{equation*}
We call {\it $q$-words} the orbits of this action.
Let $\bar w$ denote the $q$-word with representative $w$.
We call a positive integer $k$ a {\it period}
if there is an element $w_0$
in $(\mathbb Z/q\mathbb Z \times X)^*$
of length $k$ and $r_i$ in $\mathbb Z/q\mathbb Z$ for $1\le i \le |w|/k$,
such that
$$w=(r_1\cdot w_0)\cdots (r_{|w|/k} \cdot w_0).$$
Then aperiodic words are defined as usual.
That is,
$w$ is aperiodic if it is a word whose period equals its length.
With this notation, we define
$$[w]=\prod_{1\le i\le n} a_i$$
for $w=\left((r_1,a_1),\cdots,(r_n,a_n)\right).$
Then it can be easily verified that
\begin{equation}\label{definition of necklace ring*}
M^q_{\hat C}(X,\hat C^n)=\sum_{w}\,[w]\,,
\end{equation}
where the sum is over the equivalence classes of aperiodic $q$-words
$w$ out of $X$ such that the length of $w$ is $n$.
Moreover, if we specialize $x_i$ into $1$ for all $i$, then
Eq. \eqref{definition of necklace ring*} implies that
$M^q_{\hat C}(m,\hat C^n)$
represents the number of the equivalence classes of aperiodic $q$-words
$w$ out of $X$ such that the length of $w$ is equal to $n$.
Use $M^q(X,n)$ instead of $M^q_{\hat C}(X,\hat C^n)$.
Eq. \eqref{simpler form of generalized necklace} says that
\begin{equation}\label{nec poly:symmetric}
M^q(X,n)=\dfrac 1n \displaystyle\sum_{d|n} {\bf \mu^q}(n,d)\,q^{d-1}\,p_{\frac nd }(X)^d,
\end{equation}
where ${\bf \mu^q}(n,d)$ represents the
$(\hat C^n,\hat C^d)$-th entry of the matrix ${\bf \mu_{\hat C}^q}$.
More precisely,
${\bf \zeta_{\hat C}^q}$ is defined to be
\begin{equation*}
{\bf \zeta_{\hat C}^q}
(\hat C^{d_1},\hat C^{d_2})=
\begin{cases}
q^{\frac {d_1}{d_2}-1} &\text{ if } d_2 \,|\,d_1,\\
0 & \text{ otherwise},
\end{cases}
\end{equation*}
and ${\bf \mu^q_{\hat C}}$ represents its inverse.
Specializing $x_i$ into $1$ for all $i$ and then replacing $m$ by $x$,
Eq. \eqref{nec poly:symmetric} is reduced to
\begin{equation*}
M^q(x,n):=M^q_{\hat C}(x,\hat C^n)
=\dfrac 1n\sum_{d|n}{\bf \mu^q}(n,d)\,q^{d-1}\,x^d.
\end{equation*}
\section{$q$-deformation of necklace rings}\label{Necklace}
In this section, we introduce a $q$-deformation of $Nr_G$ and ${\widehat {Nr}}_G$.
Let $q$ be an indeterminate.
For a special $\ld$-ring $A$ consider the map,
\begin{equation*}
\tilde \varphi^q :A^{\mathcal O(G)} \to
{\rm Gh}(G,A),\quad
{\bf x}\mapsto \tilde\zeta_{G}^q \,{\bf x}.
\end{equation*}
For the definition of $\tilde\zeta_{G}^q$ see Eq. \eqref{def of q-matrix }.
Note that we are using the column notation.
It is obvious that if $A$ is a $\mathbb Q$-algebra, then
$\tilde\zeta_{G}^q$ is invertible and
$$(\tilde \varphi^q)^{-1}
\left(\tilde \varphi^q({\bf x})+ \tilde \varphi^q({\bf y})\right)
={\bf x}+{\bf y},$$
where the addition ${\bf x}+{\bf y}$ is defined componentwise.
\begin{lem}\label{form of product of necklace ring}
Let $A$ be a special $\ld$-ring equipped with $\mathbb Q$-algebra structure and $q$ an indeterminate.
For every $[U]\in \mathcal O(G)$, set
$$\mathfrak p^q_U:=(\tilde \varphi^q)^{-1}
\left(\tilde \varphi^q({\bf x})\cdot \tilde \varphi^q({\bf y})\right)([U])$$
Then $\mathfrak p^q_U$ is expressed as
\begin{equation*}\label{multi form of necklace ring}
\sum_{[V],[W] \preceq [U]}P_{V,W}^U(q)\,
\,\Psi^{(V:U)}({\bf x}([V]))\,\,\Psi^{(W:U)}({\bf y}([W]))
\end{equation*}
for some polynomials $P_{V,W}^U(q)\in\mathbb Q[q]$.
\end{lem}

\noindent{\bf Proof.}
For all $[Z]\in \mathcal O(G)$ we have
\begin{align*}
&\left(\tilde \zeta^q_G{\bf x}
\cdot \tilde\zeta^q_G{\bf y}\right)([Z])\\
&=\sum_{[V]\preceq [Z]\atop [W]\preceq [Z]}
\varphi_Z(G/V)\,\varphi_Z(G/W)\,q^{(V:Z)-1}q^{(W:Z)-1}\,
\Psi^{(V:Z)}({\bf x}([V])) \Psi^{(W:Z)}({\bf y}([W])).
\end{align*}
Thus, the $[U]$-th component of
$\tilde \mu^q_G\left(\tilde\zeta^q_G{\bf x}
\cdot \tilde\zeta^q_G{\bf y}\right)$
is given by
\begin{equation}\label{integraility of q-necklace}
\begin{aligned}
&\sum_{[Z]\preceq [U]}
\tilde\mu^q_G([U],[Z])\varphi_Z(G/V)\,\varphi_Z(G/W)\,q^{(V:Z)-1}q^{(W:Z)-1}
\\
&\times\left(\sum_{[V]\preceq [Z]\atop [W]\preceq [Z]}
\Psi^{(V:Z)}({\bf x}([V])) \Psi^{(W:Z)}({\bf y}([W]))\right).
\end{aligned}
\end{equation}
Exploting Lemma \ref{q-simpler form of matrix}
and the property $\Psi^m\circ \Psi^n=\Psi^{mn}, (m,n \in \mathbb N)$,
we can rewrite Eq. \eqref{integraility of q-necklace} as
\begin{align*}
&\sum_{[Z]\preceq [U] \atop {[V]\preceq [Z]\atop [W]\preceq [Z]}}
\mu^q_G([U],[Z])\varphi_Z(G/V)\,\varphi_Z(G/W)\,q^{(V:Z)-1}q^{(W:Z)-1}\,\\
&\qquad \times \,\Psi^{(V:U)}({\bf x}([V])) \Psi^{(W:U)}({\bf y}([W])).
\end{align*}
For $[V],[W]\preceq [U]$, set
$$P_{V,W}^U(q):=\sum_{[V],[W]\preceq [Z] \preceq [U]}
\mu^q_G([U],[Z])\varphi_Z(G/V)\,\varphi_Z(G/W)\,q^{(V:Z)-1}q^{(W:Z)-1}
.$$
Clearly
$$\mathfrak p^q_U=
\sum_{[V],[W] \preceq [U]}P_{V,W}^U(q)\,
\Psi^{(V:U)}({\bf x}([V]))\Psi^{(W:U)}({\bf y}([W])).$$
Since all the entries of $\mu^q_G$ are in $\mathbb Q[q]$,
we can conclude that
$$P_{V,W}^U(q)\in \mathbb Q[q].$$
\qed

\begin{lem}\label{form of product of necklace ring**}
For every $[U],[V],[W] \in \mathcal O(G)$ satisfying $[V],[W] \preceq [U]$,
$P_{V,W}^U(q)$ is a numerical polynomial in $q$.
\end{lem}
For the explicit form of $P_{V,W}^U(q)$ refer to
Example \ref{explicit of p}.
The proof of Lemma \ref{form of product of necklace ring**} will appear
in Section \ref{Proofs}.
Lemma \ref{form of product of necklace ring**} has an amusing consequence
that it yields a functor from the category of special $\ld$-rings
to the category of commutative rings.
In particular, when $q=1$, it coincides with the functor $Nr_G$
mentioned in Section \ref{preliminaries}.
\begin{thm}\label{necklace ring for special lambda rings}
Let $q$ be an integer and $G$ a profinite group.
Then there exists a unique covariant functor
$Nr_G^q$ from the category of special $\ld$-rings
to the category of commutative rings satisfying
the following conditions$:$
\begin{enumerate}
\item
As a set
$$Nr_G^q(A)=A^{\mathcal O(G)}.$$
\item
For every special $\ld$-ring homomorphism $f:A\to B$
and every ${\bf x}\in Nr_G^q(A)$ one has
$$Nr_G^q(f)({\bf x})=f\circ {\bf x}.$$
\item
The map,
\begin{equation*}
\tilde\varphi^q: Nr_G^q(A) \to {\rm Gh}(G,A),\quad
{\bf x} \mapsto \tilde \zeta_G^q \,{\bf x},
\end{equation*}
is a ring homomorphism. Here, we understand ${\bf x}$
as a $1\times \mathcal O(G)$ column vector.
\end{enumerate}
\end{thm}
Similarly, let us consider the map,
\begin{equation*}
\varphi^q :A^{\mathcal O(G)} \to
{\rm Gh}(G,A),\quad
{\bf x}\mapsto \zeta_{G}^q \,{\bf x}.
\end{equation*}
For the definition of $\zeta_{G}^q$ see Eq. \eqref{def of q-matrix: standard}.
If $A$ is a $\mathbb Q$-algebra, then
$\tilde\zeta_{G}^q$ is invertible and
$$(\varphi^q)^{-1}
\left(\varphi^q({\bf x})+ \tilde \varphi^q({\bf y})\right)
={\bf x}+{\bf y}.$$
For every $[U]\in \mathcal O(G)$, set
$$p^q_U:=(\varphi^q)^{-1}
\left(\varphi^q({\bf x})\cdot \varphi^q({\bf y})\right)([U]).$$
It is not difficult to show that $p^q_U$ is given by
\begin{equation*}
\sum_{[V],[W] \preceq [U]}P_{V,W}^U(q)\,
\,{\bf x}([V])\,\,{\bf y}([W])
\end{equation*}
for $P_{V,W}^U(q)$ appearing in Lemma \ref{form of product of necklace ring}.
Since $P_{V,W}^U(q)$'s are numerical polynomials we can state the following theorem.

\begin{thm}\label{q-deformation of necklace ring:cartier}
Let $q$ be an integer and $G$ a profinite group.
Then there exists a unique functor
${\widehat {Nr}}_G^q$ from the category of commutative rings with identity
to the category of commutative rings satisfying
the following conditions{\rm :}
\begin{enumerate}
\item
As a set
$${\widehat {Nr}}_G^q(A)=A^{\mathcal O(G)}.$$
\item
For every ring homomorphism $f:A\to B$
and every ${\bf x}\in {\widehat {Nr}}_G^q(A)$ one has
$${\widehat {Nr}}_G^q(f)({\bf x})=f\circ {\bf x}.$$
\item
The map,
\begin{equation*}
\hat \varphi^q: {\widehat {Nr}}_G^q \to {\rm Gh}(G,A),\quad
{\bf x} \mapsto \zeta_G^q \,{\bf x},
\end{equation*}
is a ring homomorphism.
\end{enumerate}
\end{thm}
As in the classical case, $Nr_G^q$ and ${\widehat {Nr}}_G^q$ will
be equivalent on the category of binomial rings.
\begin{exm}{\rm
Let $G$ be an abelian profinite group.
Then the homomorphism
$\tilde\varphi^q:Nr_G^q(A)
\to {\rm Gh}(G,A)$ is given by
\begin{equation*}
{\bf x}\mapsto
\left(\underset{ U\subseteq V} {\sum}
(G:V)\,q^{(V:U)-1}\Psi^{(V:U)}({\bf x}(V))\right)_{U\in \mathcal O(G)},
\end{equation*}
and $\varphi^q:Nr_G^q(A) \to {\rm Gh}(G,A)$ is given by
\begin{equation*}
{\bf x}\mapsto
\left(\underset{ U\subseteq V} {\sum}
(G:V)\,q^{(V:U)-1}{\bf x}(V)\right)_{U\in \mathcal O(G)}.
\end{equation*}
}\end{exm}
\begin{rem}{\rm
Recall that, in case $q=1$, ${\widehat {Nr}}_G(\mathbb Z)$
is isomorphic to the complete Burnside ring $\hat \Omega(G)$
(\cite {DS2}).
Similarly, we can realize ${\widehat {Nr}}_G^q(\mathbb Z)$
in terms of {\it twisted } Burnside ring $\hat \Omega^q(G)$.
The underlying set of $\hat \Omega^q(G)$ is same to $\hat \Omega(G)$.
Let $X$ and $Y$ be almost finite $G$-sets.
Then, from the observation
$$\tilde\varphi_U^q(X)
=\varphi_U \left((\mathbb Z/q \mathbb Z \times X) /\sim\right),$$
it follows that
\begin{align*}
&[X]\oplus[Y]=[( \mathbb Z/q\mathbb Z \times X )/\sim]
+[(\mathbb Z/q\mathbb Z \times Y)/\sim]\\
&[X]\otimes[Y]=[(\mathbb Z/q\mathbb Z \times X)/\sim]
\cdot[(\mathbb Z/q\mathbb Z \times Y)/\sim],
\end{align*}
where $\oplus,\,  \otimes$ (resp. $+,\cdot$) are operations in
${\widehat {Nr}}_G^q(\mathbb Z)$ (resp. $\hat \Omega(G)$).
Extending these operations to $\hat \Omega(G)$
we obtain a ring.
Denote this ring by $\hat \Omega^q(G)$.
By construction
$$\hat \Omega^q(G)\cong {\widehat {Nr}}_G^q(\mathbb Z).$$
}\end{rem}

\section{$q$-Inductions and $q$-restrictions on $Nr^q$}\label{on indiction and restriction}
Let $G$ be a profinite group and $U$ an open subgroup of $G$.
In this section we introduce two natural transformations
\begin{align*}
&q\text{-Ind}_U^G:Nr_U^q\to Nr_G^q,\\
&q\text{-}{\rm Res}_U^G: Nr_G^q \to  Nr_U^q,
\end{align*}
which may be viewed as a $q$-version of inductions and restrictions at $q=1$.
For the case where $q=1$ refer to \cite{O,O1}.
We start by reviewing the classical case.
References are \cite{O,O1,O2}.
\subsection{$q$-induction}
Let $A$ be a special $\ld$-ring $A$.
Then the classical induction,
$$\text{Ind}_U^G:Nr_U(A)\to Nr_G(A),\quad
{\bf x}\mapsto  \text{Ind}_U^G({\bf x}),
$$
is defined so that the $[W]$-th component of $\text{Ind}_U^G({\bf x})([W])$ equals
$$
\sum_{[V]\in \mathcal O(U)
\atop [V]=[W] \text { in } \mathcal O(G)}{\bf x}([V]).$$
Denoting  by $I_U^G$ the matrix representing $\text{Ind}_U^G$, it is immediate that
$$I_U^G
=\left(a_{[W],[V]}\right)_{[W]\in\mathcal O(G), [V]\in\mathcal O(U)},$$
where
\begin{equation*}
a_{[W],[V]}=
\begin{cases}
1& \text{ if }[V]=[W] \text{ in }\mathcal O(G),\\
0 &\text{ otherwise.}
\end{cases}
\end{equation*}
\begin{exm}{\rm
If $G$ is abelian, then
\begin{equation*}
a_{W,V}=
\begin{cases}
1& \text{ if }V=W ,\\
0 &\text{ otherwise.}
\end{cases}
\end{equation*}
In particular, when $G=\hat C$,
the operator
$$V_r:=I_{\hat C^r}^{\hat C}$$
is given by the matrix $\left(a_{i,j}\right)_{i,j\in \mathbb N}$
where
\begin{equation*}
a_{i,j}:=
\begin{cases}
1& \text{ if }(i,j) \text{ is of the form }(nr,n),\\
0 &\text{ otherwise.}
\end{cases}
\end{equation*}}
\end{exm}
\vskip 3mm
We now define $q$-induction
$$q\text{-}{\rm Ind}_U^G: Nr_U^q \to  Nr_G^q$$
in two steps. First
We will define induction
$$\nu_U^G: {\rm Gh}(U,\cdot)
\to {\rm Gh}(G,\cdot)$$
satisfying
\begin{equation}\label{how to define ind on ghost}
\tilde \varphi \circ{\rm Ind}_U^G=\nu_U^G \circ \tilde \varphi.
\end{equation}
Once $\nu_U^G$ is defined, we will show that every entry of the matrix representing
$$(\tilde \varphi^q)^{-1}\circ \nu_U^G\circ \tilde \varphi^q$$
takes its value in $\mathbb Z$.
Finally, we will define $q\text{-}{\rm Ind}_U^G$ by the multiplication by this matrix.

\vskip 3mm
{\it Step 1}:\\
Assume that $A$ is a $\mathbb Q$-algebra.
Note that the matrix representing
$${\tilde \varphi}^{-1}\circ {\rm Ind}_U^G \circ \tilde \varphi:Nr_U \to  Nr_G$$
is given by
$$\tilde\zeta_G \, I_U^G \, \tilde\mu_U
\quad \left(\stackrel{{\rm def}}{=}\left(c_{[W],[V]}\right)_{[W]\in\mathcal O(G), [V]\in\mathcal O(U)}\right).$$
Let
$$\{\varepsilon_{[V]}\,:\,[V]\in \mathcal O(U)\}
\quad \text{(resp. }\{\varepsilon_{[W]}\,:\,[W]\in \mathcal O(G)\})$$
be the standard basis of ${\rm Gh}(U,A)$
(resp. ${\rm Gh}(G,A)$), that is,
$$\varepsilon_{[Z]}([Z'])=\delta_{[Z],[Z']}.$$
Here, $\delta$ represents Kronecker's delta.
Similarly, we denote by
$$e_{[W]}\in Nr_G(A)\quad \text{(resp. } e_{[V]}\in Nr_U(A))$$
the inverse image of $\tilde \varphi$ for
$$\varepsilon_{[W]}\quad \text{(resp. }\varepsilon_{[V]}).$$

\begin{lem}\label{restriction of e(v)}{\rm (\cite[page 21]{O})}
Let $A$ be a special $\ld$-ring and also a $\mathbb Q$-algebra.
Let $U$ be an open subgroup of $G$.
Then,
$${\rm Res}_W^U(e_{[V]})=0$$
unless $[V]\preceq [W]$ in $\mathcal O(U)$.
\end{lem}

\begin{lem}\label{lemma of invariant elements}
Let $W,U$ be open subgroups of $G$ and $V$ be an open subgroup of $U$.
Consider a coset-space decomposition
$$G=\underset{i} {\dot\bigcup}\,\, g_iU$$ of $U$ in $G$, where
$\{g_i\,:\,1\le i \le (G:U)\}$ is a set of coset representatives.
If $W$ is conjugate to $V$ in $G$,
then the number of $g_i$'s satisfying the conditions,
\begin{enumerate}
\item
$g_i^{-1}Wg_i \subseteq U$
\item
{\rm (ii)}
$[g_i^{-1}Wg_i]=[V] \text{ in } \mathcal O(U)$,
\end{enumerate}
is given by
$$[N_G(V):N_U(V)].$$
\end{lem}

\noindent{\bf Proof.}
To begin with, we recall the identity
$$[N_G(V)/N_U(V)]=[N_G(V)U/U].$$
By assumption there exists an element $t\in G$ such that
$W=t^{-1}Vt$.
For our goal it suffices to show that
$$\underset{g_i^{-1}Wg_i
\stackrel{U\text{-conjugate}}{\sim}V}{\dot\bigcup} tg_iU=N_G(V)U,$$
which can be easily verified.
\qed
\vskip 3mm
In view of Lemma \ref{restriction of e(v)} and Lemma \ref{lemma of invariant elements}
we can establish the following significant result.
\begin{thm}{\rm (cf. \cite [Theorem 5.4.10]{Ben})}
\label{lem for explicit ind}
Let $A$ be a special $\ld$-ring and also a $\mathbb Q$-algebra.
Suppose that $U$ is an open subgroup of $G$ and $V$ an open subgroup of $U$.
Then
\begin{equation}
{\rm Ind}_U^G(e_{[V]})= [N_G(V):N_U(V)]e_{[V]},
\end{equation}
where the first $e_{[V]}$ is in $Nr_U(A)$ and the second in $Nr_G(A)$.
\end{thm}

\noindent{\bf Proof.}
By the Mackey formula, for an open subgroup $W$ of $G$
\begin{align*}
&\text{\rm Res}_W^G \circ \text{\rm Ind}_U^G(e_{[V]})\\
&=\sum_{WgU\subseteq G}\text{\rm Ind}_{W\cap gUg^{-1}}^W
\circ \text{\rm Res}_{W\cap gUg^{-1}}^U (g)(e_{[V]})\\
&=\sum_{WgU\subseteq G \atop [W\cap gUg^{-1}]\succeq [V] }
\text{\rm Ind}_{W\cap gUg^{-1}}^W
\circ \text{\rm Res}_{W\cap gUg^{-1}}^U (g)(e_{[V]})
\quad {\rm \text{(by Lemma \ref{restriction of e(v)})}}.
\end{align*}
Since the value of $\tilde\varphi_W$ on an element
induced from a proper open subgroup of $W$ is zero,
\begin{align*}
\tilde\varphi_W (\text{\rm Res}_W^G \circ \text{\rm Ind}_U^G(e_{[V]}))
&=\tilde\varphi_W \left(\sum_{WgU\subseteq G
\atop {W \subseteq gUg^{-1}
\atop [gUg^{-1}]\succeq [V]}}
\text{\rm Res}_{gUg^{-1}}^U (g)(e_{[V]})\right)\\
&=\begin{cases}
\displaystyle\sum_{WgU\subseteq G
\atop {g^{-1}Wg \subseteq U
\atop [g^{-1}Wg]=[V] \text{ in } \mathcal O(U)}}1
&\text { if } W\stackrel{\rm conjugate}{\sim} V\,,\\
0&\text { otherwise.}
\end{cases}\\
&=\begin{cases}
[N_G(V):N_U(V)]
&\text { if } W\stackrel{\rm conjugate}{\sim} V,\\
0&\text { otherwise.}
\end{cases}
\end{align*}
The last equality follows from Lemma \ref{lemma of invariant elements}.
Hence, we have
\begin{equation*}
{\rm Ind}_U^G(e_{[V]})= [N_G(V):N_U(V)]e_{[V]}.
\end{equation*}
\qed
\vskip 3mm
Theorem \ref{lem for explicit ind} implies that
$(c_{[W],[V]})$ is given by
\begin{equation}\label{explicit matrix form**}
c_{[W],[V]}=
\begin{cases}
[N_G(V):N_U(V)]& \text{ if }[V]=[W] \text{ in }\mathcal O(G),\\
0 &\text{ otherwise.}
\end{cases}
\end{equation}
Now, for an arbitrary special $\ld$-ring $A$, define
$$\nu_U^G: {\rm Gh}(U,A)
\to {\rm Gh}(G,A)$$
by
$${\bf x}\mapsto (c_{[W],[V]}){\bf x}$$
for all ${\bf x} \in {\rm Gh}(U,A).$
By definition it is straightforward that
$$\tilde \varphi\circ {\rm Ind}_U^G
=\nu_U^G\circ \tilde \varphi.$$
\begin{exm}{\rm
If $G$ is abelian, then
\begin{equation*}
c_{V,W}=
\begin{cases}
(G:U)& \text{ if }W=V ,\\
0 &\text{ otherwise.}
\end{cases}
\end{equation*}
In particular, when $G=\hat C$,
the operator $V_r:=\nu_{\hat C^r}^{\hat C}$
is given by the matrix $\left(a_{i,j}\right)_{i,j\in \mathbb N}$
where
\begin{equation*}
c_{i,j}:=
\begin{cases}
r& \text{ if }(i,j) \text{ is of the form }(nr,n),\\
0 &\text{ otherwise.}
\end{cases}
\end{equation*}
}\end{exm}

\vskip 3mm
{\it Step 2}:\\
Now, let us define $q$-induction
$$q\text{-}{\rm Ind}_U^G: Nr_U^q(A) \to  Nr_G^q(A).$$
In case where $A$ is a $\mathbb Q$-algebra,
it can be defined by
\begin{equation*}
(\tilde \varphi^q)^{-1}\circ \nu_U^G \circ \tilde \varphi^q.
\end{equation*}
It is quite interesting to see that
$q$-induction thus defined coincides with the classical one
for every open subgroup $U$ of $G$. We need the following lemma.
\begin{lem}\label{two same formulas}\hfill

{\rm (a)}
For open subgroups $W,U$ of $G$ we have
$$\varphi_W(G/U)=[N_G(U):U]\,n(W,U)
=\sum_{[W_i]\in \mathcal O(U)\atop [W_i] =[W]\,\,\, \text{\rm in
}\mathcal O(G)}[N_G(W_i):N_U(W_i)],$$
where $n(W,U)$ is the number
of $G$-conjugates of $U$ containing $W$.

{\rm (b)}
For open subgroups $W,U$ of $G$
and an open subgroup $V$ of $U$,
we have
$$\varphi_W(G/V)
=\sum_{[W_i]\in \mathcal O(U)\atop [W_i]
=[W]\,\,\, \text{\rm in }\mathcal O(G)}[N_G(W_i):N_U(W_i)]\varphi_{W_i}(U/V).$$
\end{lem}

\noindent{\bf Proof.}
(a) In fact the first equality is well-known in case $G$ is finite
(for example, see \cite {D}), and it is easy to check that
it is also true in case $G$ is
an arbitrary profinite group.
In detail, let $\{g_i: 1\le i \le (G:U)\}$ be a complete set of
left coset representatives of $U$ in $G$.
Then $\varphi_W(G/U)$ is the number of $g_i$'s satisfying
\begin{equation*}\label{first equation}
W\subseteq g_i U g_i^{-1}.
\end{equation*}
Note that $g_i\in g_j N_G(U)$ if and only if $g_i U g_i^{-1}=g_j U g_j^{-1}$.
Hence, $[N_G(U):U]$-number of $g_i$'s yields one $G$-conjugate.
From this observation the first equality follows.

The second equality follows from \eqref{explicit matrix form**}.
By definition of the induction map and $\tilde \varphi_W$ we have
\begin{equation*}
\tilde \varphi\circ {\rm Ind}_U^G(\varepsilon_{[U]})
=\varphi_W(G/U).
\end{equation*}
On the other hand, \eqref{explicit matrix form**} yields that
\begin{equation*}
\nu_U^G\,( \tilde \varphi(\varepsilon_{[U]}))([W])
=\sum_{[W_i]\in \mathcal O(U)\atop [W_i]=[W]\,\,\,
\text{\rm in }\mathcal O(G)}[N_G(W_i):N_U(W_i)].
\end{equation*}
So we are done.

(b) can be proved in the exactly same way as in (a).
\qed

\begin{thm}\label{on induction}
Let $A$ be a special $\ld$-ring equipped with $\mathbb Q$-algebra structure.
Then,
$$q \text{-}{\rm Ind}_U^G={\rm Ind}_U^G$$
for every open subgroup $U$ of $G$ we have
\end{thm}

\noindent{\bf Proof.}
It suffices to show
\begin{equation*}
\tilde \varphi^q\circ {\rm Ind}_U^G(\varepsilon_{[V]})
=\nu_U^G \circ \tilde \varphi^q(\varepsilon_{[V]})
\end{equation*}
for every $[V]\in \mathcal O(U)$.
Note that for $[W]\in \mathcal O(G)$
\begin{align*}
&\tilde \varphi^q\circ {\rm Ind}_U^G(\varepsilon_{[V]})([W])\\
&=\varphi_W(G/V)\,q^{(V:W)-1}\\
&=\sum_{[W_i]\in \mathcal O(U)\atop [W_i]
=[W]\,\,\, \text{\rm in }\mathcal O(G)}[N_G(W_i):N_U(W_i)]
\varphi_{W_i}(U/V)q^{(V:W)-1}\text{ (by Lemma } \ref{two same formulas})\\
&=\nu_U^G \circ \tilde \varphi^q(\varepsilon_{[V]})([W]).
\end{align*}
So we are done.
\qed
\vskip 3mm
Theorem \ref{on induction} enables us to generalize $q$-induction for arbitrary special $\ld$-rings
in the obvious way.
More precisely,
$$q\text{-}{\rm Ind}_U^G:
 Nr_U^q(A) \to  Nr_G^q(A),
\quad {\bf x} \mapsto I_U^G\,{\bf x},$$
where $I_U^G$ is the matrix representing $\text{Ind}_U^G$.
By definition it is obvious that
\begin{equation*}
\tilde \varphi^q\circ q\text{-}{\rm Ind}_U^G
=\nu_U^G\circ \tilde \varphi^q.
\end{equation*}
Similarly, if we define
$$q\text{-Ind}_U^G: {\widehat {Nr}}_U^q \to  {\widehat {Nr}}_G^q$$
by the multiplication by $I_U^G$, then
\begin{equation*}
\hat  \varphi^q\circ q\text{-}{\rm Ind}_U^G
=\nu_U^G\circ \hat \varphi^q.
\end{equation*}

\subsection{$q$-restriction}\label{q-restriction: subsection}
In this section, we investigate how to define $q$-restriction.
In contrast with $q$-induction,
$q$-restriction depends on $q$.
First, we recall the classical case, i.e., the case
where $q=1$. References are \cite{O,O1}.
Given a special $\ld$-ring $A$,
restriction,
$$\text{\rm Res}_U^G : Nr_G(A)\to Nr_U(A),
\quad {\bf x} \mapsto \text{\rm Res}_U^G({\bf x}),
$$
is defined by the assignment
$$\text{\rm Res}_U^G({\bf x})([W])
=\underset{[V]\in \mathcal O(G)}{{\sum}}\sum_{[Z(g,U,V)]
=[W]}\Psi^{(V\,:\,Z(g,U,V))}({\bf x}([V])).$$
Here, $g$ ranges over the set
of representatives of $U$-orbits of $G/V$
and $Z(g,U,V)$ stands for $U\cap gVg^{-1}$.
Also, restriction on ghost rings is defined as follows:
\begin{equation*}
\mathcal F_U^G: {\rm Gh}(G,A)\to{\rm Gh}(U,A),
\quad
{\bf x} \mapsto R_U^G\cdot {\bf x},
\end{equation*}
where $R_U^G$ represents the $\mathcal O(U)\times \mathcal O(G)$
matrix given by
$$\left(b_{[V],[W]}\right)_{[V]\in\mathcal O(G), [W]\in\mathcal O(U)}$$
with
\begin{equation*}
b_{[V],[W]}=
\begin{cases}
1& \text{ if }[W]=[V] \text{ in }\mathcal O(G),\\
0 &\text{ otherwise.}
\end{cases}
\end{equation*}
\begin{lem}{\rm (\cite{O})}
For any special $\ld$-ring $A$,
$$\tilde \varphi^q\circ \text{\rm Res}_U^G=\mathcal F_U^G\circ\tilde \varphi^q.$$
\end{lem}

\begin{exm}{\rm
In case where $G$ is abelian,
$R_U^G$ is given by the matrix
$\left(b_{V,W}\right)_{V,W }$
with
\begin{equation*}
b_{V,W}=
\begin{cases}
1& \text{ if }V=W,\\
0 &\text{ otherwise.}
\end{cases}
\end{equation*}
In particular, if $G=\hat C,$
\begin{align*}
&\mathcal F_{\hat C^r}^{\hat C}:A^{\mathbb N} \to A^{\mathbb N}\\
&\qquad {\bf x} \to R_r {\bf x},
\end{align*}
where $R_r(:=R_{\hat C}^{\hat C^r})$
is given by the matrix $\left(a_{i,j}\right)_{i,j\in \mathbb N}$
with
\begin{equation*}
a_{i,j}=
\begin{cases}
1& \text{ if }(i,j) \text{ is of the form }(n,nr),\\
0 &\text{ otherwise.}
\end{cases}
\end{equation*}
}\end{exm}

\vskip 3mm
Let us first assume that
$A$ is a $\mathbb Q$-algebra.
Set
\begin{equation*}
q\text{-}{\rm Res}_U^G=(\tilde \varphi^q)^{-1}
\circ \mathcal F_U^G\circ \tilde \varphi^q.
\end{equation*}

\begin{lem}\label{Mackey subgroup theorem}
Let $A$ be a special $\ld$-ring
equipped with $\mathbb Q$-algebra structure.
For any two open subgroups $U,\,V $ of $G$,
the composite map
$$q\text{\rm -Res}_U^G \circ \text{\rm Ind}_V^G :
Nr^q_V(A)\to Nr^q_U(A)$$ is given by
\begin{equation}\label{q:Mackey subgroup theorem}
q\text{\rm -Res}_U^G \circ \text{\rm Ind}_V^G
=\sum_{UgV\subseteq G}\text{\rm Ind}_{U\cap gVg^{-1}}^U
\circ
q\text{\rm -Res}_{U\cap gVg^{-1}}^V (g).
\end{equation}
\end{lem}

\noindent{\bf Proof.}
The left-hand side of
Eq. \eqref{q:Mackey subgroup theorem} equals
\begin{align*}
q\text{\rm -Res}_U^G \circ \text{\rm Ind}_V^G
&=(\tilde \varphi^q)^{-1}
\circ \mathcal F_U^G\circ \tilde \varphi^q \circ \text{\rm Ind}_V^G\\
&=(\tilde \varphi^q)^{-1}
\circ \mathcal F_U^G\circ \nu_U^G \circ \tilde \varphi^q,
\end{align*}
and the right-hand side of Eq. \eqref{q:Mackey subgroup theorem}
equals
\begin{align*}
&\sum_{UgV\subseteq G}\text{\rm Ind}_{U\cap gVg^{-1}}^U
\circ q\text{\rm -Res}_{U\cap gVg^{-1}}^V (g)\\
&=\sum_{UgV\subseteq G}\text{\rm Ind}_{U\cap gVg^{-1}}^U
\circ (\tilde \varphi^q)^{-1}
\circ \mathcal F_{U\cap gVg^{-1}}^V(g)\circ \tilde \varphi^q\\
&=\sum_{UgV\subseteq G}(\tilde \varphi^q)^{-1} \circ
\nu_{U\cap gVg^{-1}}^U
\circ \mathcal F_{U\cap gVg^{-1}}^V(g)\circ \tilde \varphi^q.
\end{align*}
Thus, for our purpose, it suffices to show that
\begin{equation*}
\mathcal F_U^G\circ \nu_U^G
=\sum_{UgV\subseteq G}\nu_{U\cap gVg^{-1}}^U
\circ \mathcal F_{U\cap gVg^{-1}}^V(g).
\end{equation*}
This identity is also equivalent to
\begin{equation*}
\text{\rm Res}_U^G \circ \text{\rm Ind}_V^G
=\sum_{UgV\subseteq G}\text{\rm Ind}_{U\cap gVg^{-1}}^U
\circ \text{\rm Res}_{U\cap gVg^{-1}}^V (g),
\end{equation*}
which follows from \cite {O}.
\qed
\vskip 3mm

Consider the matrix representing the map $q\text{-}{\rm Res}_U^G$
\begin{equation*}
\tilde\mu_{U}^q \, R_U^G \, \tilde\zeta_{U}^q
\quad (\stackrel{{\rm def}}{=}(d_{[V],[W]}(q))_{[V]\in\mathcal O(U), [W]\in\mathcal O(G)}).
\end{equation*}
By direct computation we obtain that
\begin{align*}
d_{[V],[W]}(q)=
\begin{cases}
Q_{V,W}(q)\Psi^{(W:V)}&\text{ if } [W]\preceq [V] \text{ in }
\mathcal O(G),\\
0&\text{ otherwise,}
\end{cases}
\end{align*}
where
\begin{equation}\label{q-rest}
Q_{V,W}(q)=\sum_{[S]\in \mathcal O(U)
\atop {[W]\preceq [S]\text{ in }\mathcal O(G)
\atop  [S]\preceq
[V]\text{ in }\mathcal O(U)}}
\mu_U^q([V],[S])\varphi_S(G/W)q^{(W:S)-1}.
\end{equation}

\begin{thm}\label{well of q-restric}
Let $A$ be a special $\ld$-ring with $\mathbb Q$-algebra structure. Then,
$Q_{V,W}(q)$ is a numerical polynomial in $q$ for every
$[V]\in\mathcal O(U)$ and every $[W]\in\mathcal O(G).$
\end{thm}
The proof of Theorem \ref{well of q-restric} will appear in Section \ref{proof of q-rest}.
As a direct consequence of Theorem
\ref{well of q-restric} we will extend $q$-restriction to arbitrary special $\ld$-ring
in the obvious way. More precisely, given a special $\ld$-ring $A$, we
define
$$q\text{-}{\rm Res}_U^G:
 Nr_G^q(A)\to Nr_U^q(A)$$
by
$${\bf x} \mapsto(d_{[V],[W]}(q))({\bf x}).$$
By definition it is obvious that
$$\tilde \varphi^q\circ q\text{-}{\rm Res}_U^G
=\mathcal F_U^G \circ \tilde \varphi^q.$$
Similarly, if we define
$$q\text{-Res}_U^G: {\widehat {Nr}}_G^q \to  {\widehat {Nr}}_U^q$$
by the multiplication by the matrix
$$\left(Q_{V,W}(q)\right)_{[V]\in\mathcal O(U), [W]\in\mathcal O(G)},$$
then
$$\hat \varphi^q\circ q\text{-}{\rm Res}_U^G
=\mathcal F_U^G \circ \hat \varphi^q.$$
\section{$q$-deformation of Witt-Burnside rings}\label{Witt-Burn}
In this section, we construct a $q$-deformation
of the Witt-Burnside ring of a profinite group, where $q$ ranges over the set of integers.
Let $G$ be a profinite group. Put
$$A=\mathbb Q\,[\,{\bf x}([U]),{\bf y}([U]): [U] \in \mathcal O(G)\,],$$
where ${\bf x}([U]),{\bf y}([U])$ are indeterminates.
Consider the map
\begin{align*}
&\Phi^q:A^{\mathcal O(G)}\to {\rm Gh}(G,A),\\
&{\bf x}\mapsto
\left( \underset{[V] \preceq [U]}{{\sum}}
\varphi_U(G/V)\,q^{(V:U)-1}\, {\bf x}([V])^{(V:U)}\right)_{[U]\in \mathcal O(G)}.
\end{align*}
For every $[U]\in \mathcal O(G)$ we let
$$\mathfrak s_{[U]}:=(\Phi^q)^{-1}(\Phi^q({\bf x})+\Phi^q({\bf y}))([U])$$
and
$$\mathfrak p_{[U]}:=(\Phi^q)^{-1}(\Phi^q({\bf x})\cdot \Phi^q({\bf y}))([U]).$$
With this notation we have

\begin{lem}\label{q-deformation of witt-Burn}
Let $U$ be an open subgroup of $G$.
Under the above hypothesis, we have
$$\mathfrak s_{[U]},\,\mathfrak p_{[U]}\in \mathbb Z\,[\,{\bf x}([V]),{\bf y}([V])\,:\,[V] \preceq [U]\,].$$
Furthermore, given ${\bf x}\in A^{\mathcal O(G)}$, there exists a unique
$\iota_{{\bf x}}$ such that
$$\iota_{{\bf x}}([U])\in \mathbb Z\,[\,{\bf x}([V]):\,[V] \preceq [U]\,]$$
and
$$0=s_{[U]}({\bf x}([V]),\iota_{{\bf x}}([V]) \,:\, [V]\preceq [U]).$$
\end{lem}

\begin{exm}
{\rm
Let $G$ be an abelian profinite group
and $U$ an open subgroup of $G$ with $(G:U)=p$, where $p$ is a prime.
Clearly $\mathfrak s_{G}={\bf x}(G)+{\bf y}(G)$ and
$\mathfrak p_{G}={\bf x}(G){\bf y}(G).$
Thus,
\begin{align*}
&\mathfrak s_{[U]}={\bf x}(U)+{\bf y}(U)-\frac {q^{p-1}}{p}\sum_{r=1}^{p-1}{p\choose r}{\bf x}(G)^r{\bf y}(G)^{p-r},\\
&\mathfrak p_{[U]}=\frac {q^{p-1}(q^{p-1}-1)}{p}{\bf x}(G)^p{\bf y}(G)^p
+q^{p-1}({\bf x}(G)^p{\bf y}(U)+{\bf x}(G){\bf y}(G)^p)+p\,{\bf x}(U){\bf y}(U).
\end{align*}
}\end{exm}

The proof of Lemma \ref{q-deformation of witt-Burn} will appear in Section \ref{Proofs}.
In view of Lemma \ref{q-deformation of witt-Burn}
one can derive a functor from the category of commutative rings
with identity to the category of commutative rings.
\begin{thm}
Let $q$ be an integer and $G$ be a profinite group.
Then there exists a unique covariant functor
${\mathbb W}_G^q$ from the category of commutative rings with identity
into the category of commutative rings satisfying
the following conditions{\rm :}
\begin{enumerate}
\item
As a set
$${\mathbb W}_G^q(A)=A^{\mathcal O(G)}.$$
\item
For every ring homomorphism $f:A\to B$ and every $\alpha
\in  \mathbb W_G^q(A)$ one has
$$\mathbb W_G^q(f)(\alpha)=f\circ \alpha.$$
\item
The map,
\begin{align*}
&\Phi^q:  \mathbb W_G^q(A) \to {\rm Gh}(G,A),\\
&\alpha \mapsto \left(\underset{[G]\preceq [V] \preceq[ U]}{\sum}
\varphi_U(G/V)\, q^{(V:U)-1}\alpha([V])^{(V:U)}\right)_{[U]\in\mathcal O(G)},
\end{align*}
is a ring homomorphism.
\end{enumerate}
\end{thm}
From the third condition it follows that
\begin{equation}\label{relatio of phi and q-phi}
\Phi^q(\alpha)=\frac 1q \Phi(q\alpha)\quad \text{ if }q\ne 0.
\end{equation}
Here, $q\alpha$ denotes the vector whose $[U]$-th component is $q\alpha([U])$.
This identity remains effective even when $q=0$
since $$\Phi(q\alpha)([U])$$ is divided by $q$ for every $[U]\in \mathcal O(G)$.

\begin{exm}{\rm
If $G$ is abelian, then
$\Phi^q: \mathbb W_G^q(A)\to  {\rm Gh}(G,A)$
is given by
\begin{equation*}
{\bf x}\mapsto \left(\underset{U\subseteq V} {\sum}
(G:V)\,q^{(V:U)-1}{\bf x}(V)^{(V:U)}\right).
\end{equation*}
In particular, if $G=\hat C$, then
$\Phi^q:
\mathbb W_{\hat C}^q(A)\to A$
is given by
\begin{align*}
{\bf x}\mapsto
\left(\sum_{d|n}\,d\,q^{\frac nd-1}{\bf x}(d)^{\frac nd}
\right).
\end{align*}
}\end{exm}

Denote by ${\rm D}(G)$ the set $\{(N_G(U):U) :[U]\in \mathcal
O(G)\}$ and by $\mathbb Z_{G}$ the commutative ring
\begin{equation}\label{defn of local}
\mathbb Z\left[\frac {1}{(N_G(U):U)}
:[U]\in \mathcal O(G)\right].
\end{equation}
It is not difficult to show the following fact.

\begin{lem}\label{injectivity and bijectivity}
Let $G$ be a profinite group.
Then we have the following characterizations.

$($a$)$
$\Phi_G^q(A)$ is injective
$\Leftrightarrow $
$\varphi_G^q(A)$ is injective
$\Leftrightarrow $
$A$ has no $\{(N_G(U):U)
:[U]\in \mathcal O(G)\}$-torsion.

$($b$)$
$\Phi_G^q(A)$ is surjective
$\Leftrightarrow $
$\varphi_G^q(A)$ is surjective
$\Leftrightarrow $
$\Phi_G^q(A)$ is bijective
$\Leftrightarrow $
$\varphi_G^q(A)$ is bijective
$\Leftrightarrow $
$A$ is a $\mathbb Z_{G}$-algebra.
\end{lem}
In the classical case we have constructed an isomorphism
$${\tau}: \mathbb W_G(A)\to Nr_G(A),$$
called {\it Teichm\"{u}ller map},
for every special $\ld$-ring $A$ (\cite{O,O2}).
In the following, we will introduce its $q$-deformation
$${\tau^q}: \mathbb W_G^q(A)\to Nr_G^q(A).$$
To do this, let us first define $q$-exponential maps
$$\tau^G_q: A\to  Nr_G^q(A).$$
Given an element $r \in A$, write
it as a sum of one-dimensional elements, say
$r_1+r_2+\cdots+r_m.$
Then, from Eq.  \eqref{cloased form: q-necklacr poly} we have
\begin{equation}\label{well-defined of q-neck}
M_G^q(\{r_1,r_2, \cdots, r_m\},V)=\sum_{[W]\preceq [V]}\mu^q_{G}(V,W)\,
\,q^{(G:W)-1}\,\Psi^{(G:W)}(r^{(W:V)}).
\end{equation}
for every open subgroup $V$ of $G$.
Set
$$M_G^q(r,V)=M_G^q(\{r_1,r_2, \cdots, r_m\},V).$$
From Eq. \eqref{well-defined of q-neck} it follows that it is well-defined, that is,
it does not depend on the choice of decompositions of $r$ into one-dimensional elements.
\begin{lem}\label{well-defineness of expo}
Let $A$ be a special $\ld$-ring.
Then, for every $r\in A$,
${M}_G^q(r,V)\in A$.
\end{lem}
\noindent{\bf Proof.}
By definition $M_G^q(r,V)$ is a
symmetric polynomial in $r_i$'s.
Hence, it is an integral
polynomial in $\ld^k(r),\,\,1\le k \le (G:V),$
which is well known in the context of symmetric functions.
Since $\ld^k(r) \in A$ for all $r\in A$ and $k\ge 1$, we have the desired result.
\qed
\vskip 3mm
Set
$$\tau^G_q(r)([V])=M_G^q(r,V)$$
for all $[V]\in \mathcal O(G).$
Lemma \ref{well-defineness of expo} implies that $\tau^G_q(r)$ is an element of
$Nr_G^q(A)$ for all $r\in A$.
Moreover,
\begin{equation}\label{comm of varphi and tau}
\tilde\varphi^q\circ \tau^G_q(r)([V])=q^{(G:V)-1}r^{(G:V)}
\end{equation}
since
$\tilde \varphi^q$ represents the multiplication by $\tilde \zeta_G^q$
(see Theorem \ref{necklace ring for special lambda rings}).
Using $q\text{-Ind}_U^{G}$
and $\tau^U_q$ for all $[U]\in \mathcal O(G)$
simultaneously, we construct
a $q$-analog of Teichm\"{u}ller map as follows:
\begin{align*}
&\tau^q:\mathbb W_G^q(A)\to Nr_G^q(A),\\
&\alpha\mapsto {\underset{[U]\in \mathcal O(G)}{\sum}}\,\text{Ind}_U^{G}
\circ \tau^U_q(\alpha([U])).
\end{align*}
We claim that $\Phi^q=\tilde\varphi^q \circ \tau^q$.
In proving this argument, the following proposition plays a crucial role.
\begin{prop}
Let $A$ be a special $\ld$-ring with $\mathbb Q$-algebra
structure.
Then, for every $[U]\in \mathcal O(G)$ and $r\in A$, the identity
\begin{equation}\label{relation of q-restriction}
q\text{-}{\rm Res}_U^G\left(\tau_q^G(r) \right)
=\tau_q^U(q^{(G:U)-1}r^{(G:U)})
\end{equation}
holds.
\end{prop}
\noindent{\bf Proof.}
From Eq. \eqref{comm of varphi and tau}
and the definition of ${\mathcal F}_U^G$ (see Section \ref{q-restriction: subsection}),
the $[W]$-th component of ${\mathcal F}_U^G\circ\tilde\varphi^q\circ \tau^G_q(r)$ is given by
\begin{align*}
\begin{cases}
q^{(G:V)-1}r^{(G:V)}& \text{ if }[W]=[V] \text{ in } \mathcal O(G) ,\\
0&\text{ otherwise,}
\end{cases}
\end{align*}
for all $[W]\in \mathcal O(G).$
Hence,
\begin{equation*}
\begin{aligned}
&q\text{-}{\rm Res}_U^{G}\left(\tau_q^U(r) \right)([W])\\
&=(\tilde \varphi^q)^{-1}\circ {\mathcal F}_U^G
\circ\tilde \varphi^q\left(\tau_q^U(r) \right)([W])\\
&=\sum_{[U]\preceq [Z]\preceq [W]}\mu_U^q([W],[Z])\,\,
q^{(G:Z)-1}\Psi^{(Z:W)}(r^{(G:Z)})\\
&=\sum_{[U]\preceq [Z]\preceq [W]}\mu_U^q([W],[Z])\,\,
q^{(U:Z)-1}\Psi^{(Z:W)}
\left(\left(q^{(G:U)-1)}r^{(G:U)}\right)^{(U:Z)}\right)\\
&=M^q_U\left(q^{(G:U)-1)}r^{(G:U)},W\right).
\end{aligned}
\end{equation*}
This completes the proof.
\qed

\begin{thm} \label{ringhomo of q-teichmuller*}
$\tau^q$ is bijective for every special $\ld$-ring $A$.
\end{thm}

\noindent{\bf Proof.}
For $\alpha,\beta \in \mathbb W_G^q(A)$ we assume that
$${\underset{[U] \in \mathcal O(G)}{\sum}}\,\text{Ind}_U^{G}
\circ \tau^U_q(\alpha([U]))
={\underset{[U] \in \mathcal O(G)}{\sum}}\,\text{Ind}_U^G
\circ \tau^U_q(\beta([U])).$$
Then, for every $[U] \in \mathcal O(G),$
\begin{equation}\label{inj of q-teich}
\left(\underset {[V] \in \mathcal O(G)\atop [U] \in \mathcal O(V)}
{{\sum}}M_V^q(\alpha([V]),U)\right)
=\left(\underset {[V] \in \mathcal O(G)\atop [U] \in \mathcal O(V)}
{{\sum}}M_V^q(\beta([V]),U)\right).
\end{equation}
It is clear that if $V=G$, then $\alpha([G])=\beta([G]).$
Now assume that $\alpha([V])=\beta([V])$ for all $V$ such that $(G:V)<(G:U).$
From \eqref{inj of q-teich} it follows that $\alpha([U])=\beta([U]).$
Thus, $\alpha=\beta,$ and which implies the injectiveness of $\tau^q$.

Next, we will show that $\tau^{q}$ is surjective.
For any ${\bf a}\in Nr_G^q(A)$
we would like to find an element ${\bf x}\in  \mathbb W_G^q(A)$
satisfying
\begin{equation}\label{surjective of teich:}
\left(\sum_{[V] \in \mathcal O(G)\atop [U] \in \mathcal O(V)}
M_V^q({\bf x}([V]),U)\right)
={\bf a}([U])
\end{equation}
for every $[U]\in \mathcal O(G)$.
If $U=G$, then ${\bf x}([G])={\bf a}([G])$.
Let us use mathematical induction on the index.
Assume that we have found ${\bf x}([V])$ for all $[V]\in \mathcal O(G)$
such that $(G:V)<(G:U).$
Setting
$$
{\bf x}([U])
={\bf a}([U])-
\left(\sum_{[V]\in \mathcal O(G),\,[V]\ne [U]\atop [U] \in \mathcal O(V)}
M_V^q({\bf x}([V]),U)\right),
$$
we have
$$
\left(\sum_{[V]\in \mathcal O(G)\atop [U] \in \mathcal O(V)}
M_V^q({\bf x}([V]),U)\right)={\bf a}([U])$$
since
$M_U^q({\bf x}([U]),U)={\bf x}([U])$.
In this way we can find ${\bf x}$ satisfying Eq. \eqref{surjective of teich:}.
\qed
\vskip 3mm
If $U$ and $V$ are open subgroups of $G$
and if $g\in G$ conjugates
$U$ into $V$, we can induce $q$-restriction and $q$-induction
$q\text{-Res}_U^V(g):Nr^q_V(A)\to Nr^q_U(A)$
and $\text{\rm Ind}_U^V(g): Nr^q_U(A)\to Nr^q_V(A)$
from the embedding $U\hookrightarrow V$.
With this notation, we can extend some classical facts on $q$-induction and $q$-restriction,
which can be proved by exploiting the facts in \cite[Section 2.11]{DS2}.
\begin{prop}\label{Mackey}\hfill

{\rm (a)} {\rm (Frobenius reciprocity) }
For any two open subgroups $U,\,V $ of $G$ and
${\bf x}\in Nr^q_G(A)$ and ${\bf y}\in Nr^q_U(A)$,
one has
$$\text{\rm Ind}_U^G({\bf y})\cdot {\bf x}
=\text{\rm Ind}_U^G({\bf y}\cdot q\text{\rm -Res}_U^G({\bf x})).$$

{\rm (b)} {\rm (Mackey subgroup theorem) }
For any two open subgroups $U,\,V $ of $G$,
the composite map
$$q\text{\rm -Res}_U^G \circ \text{\rm Ind}_V^G :
Nr^q_V(A)\to Nr^q_U(A)$$ is given by
\begin{equation}\label{q:Mackey subgroup theorem*}
q\text{\rm -Res}_U^G \circ \text{\rm Ind}_V^G
=\sum_{UgV\subseteq G}\text{\rm Ind}_{U\cap gVg^{-1}}^U
\circ
q\text{\rm -Res}_{U\cap gVg^{-1}}^V (g).
\end{equation}

{\rm (c)}
For any two open subgroups $U,\,V $ of $G$,
\begin{equation}
\label{property 3}
\tilde\varphi^q_U\circ \text{\rm Ind}_V^G =
\sum_{gV\in (G/V)^U}\tilde\varphi^q_U \circ  q\text{\rm -Res}_U^V (g).
\end{equation}
\end{prop}
\noindent{\bf Proof. }
The proofs can be done with small modification of those of \cite[Section 2.11]{DS2}.
So we will prove only (c).
Note that $\tilde \varphi^q_U=\tilde \varphi^q_U\circ q\text{\rm -Res}_U^G$
since
\begin{align*}
\tilde \varphi^q_U\circ q\text{\rm -Res}_U^G ({\bf x})
&=\mathcal F_U^G\circ \tilde \varphi^q ({\bf x})([U])\\
&=\tilde \varphi^q ({\bf x})([U]) \quad \text{ (by def of }\mathcal F_U^G).
\end{align*}
Furthermore $\tilde \varphi^q_U\circ \text{\rm Ind}_{U\cap gVg^{-1}}^U=0$
unless $U=U\cap gVg^{-1}$, equivalently, $$gV\in (G/V)^U.$$
Applying these two facts to Eq. \eqref{q:Mackey subgroup theorem*}
yields the desired result.
\qed

\begin{thm}\label{ringhomo of q-teichmuller}
The following diagram
\begin{equation*}
\begin{picture}(360,70)
\put(100,60){$\mathbb W_G^q(A)$}
\put(135,63){\vector(1,0){70}}
\put(210,60){$Nr_G^q(A)$}
\put(155,0){${\rm Gh(G,A)}$}
\put(215,55){\vector(-1,-1){43}}
\put(120,55){\vector(1,-1) {43}}
\put(165,67){$\tau^{q}$}
\put(205,35){$\tilde\varphi^q$}
\put(121,34){$\Phi^q$}
\put(163,38){$\curvearrowright$}
\end{picture}
\end{equation*}
is commutative.
\end{thm}

\noindent{\bf Proof.}
For $r\in A$,
\begin{align*}
\tilde\varphi_U^{q}
\left(\text{Ind}_V^{G}
\circ \tau^V_q(r)\right)
&=\sum_{gV\in (G/V)^U}
\tilde\varphi^q_U(q\text{-Res}_U^{V}(g)(\tau^V_q(r))
\qquad \text{ by }\eqref {property 3}\\
&=\sum_{gV\in (G/V)^U}
\tilde\varphi^q_U \left(\tau^U_q\left(q^{(V:U)-1}r^{(V:U)}\right)\right)
\qquad \text{ by } \eqref{relation of q-restriction}\\
&=\varphi_U(G/V)\,q^{(V:U)-1}r^{(V:U)}.
\end{align*}
Therefore, the additivity of $\tilde \varphi$ justifies our assertion.
\qed

\begin{thm}
For every special $\ld$-ring $A$,
the map
$\tau^{q}:\mathbb W_G^q(A)\to Nr_G^q(A)$ is a ring isomorphism.
\end{thm}
\noindent{\bf Proof.}
We have already proved the bijectiveness of $\tilde \varphi$ in Theorem \ref{ringhomo of q-teichmuller*}.
In order to show that $\tilde \varphi$ is a ring homomorphism
let $A=\mathbb Q[{\bf x}([U]),{\bf y}([U]): [U]\in \mathbb O(G)]$.
Consider
\begin{equation}\label{universal identity}
\begin{aligned}
&\tau^{q}({\bf x}+{\bf y})
=\tau^{q}({\bf x})+\tau^{q}({\bf y}).
\end{aligned}
\end{equation}
We claim that the above identity is {\it universal} in the sense that, for every $[U]\in \mathcal O(G)$,
the $[U]$-th component of both sides equals as a polynomial with integral coefficients.
More precisely,
Eq. \eqref{universal identity} gives rise to the equality
\begin{equation}\label{poly identity}
\sum_{[V]\in \mathcal O(G)}
M_G^{q}(({\bf x}+{\bf y})([V]), U)
=\sum_{[V]\in \mathcal O(G)}
\left(M_G^{q}({\bf x}([V]),U)+M_G^{q}({\bf y}([V]),U)\right)
\end{equation}
for every $[U]\in \mathcal O(G)$.
First, note that
$$M_G^q({\bf x}([V]),U),\quad [V],[U]\in
\mathcal O(G),$$
is a polynomial in $\ld^n({\bf x}([V])), n\ge 1$
with integral coefficients. Here, $\ld^n$ represents the $n$-th $\ld$-operation.
Second, observe that
$$({\bf x}+{\bf y})([V]),\quad ({\bf x}\cdot{\bf y})([V])
\in \mathbb Z\,[\,{\bf x}([W]),{\bf y}([W])\,:\,[W] \preceq [V]\,].$$
(see Lemma \ref{q-deformation of witt-Burn})
Putting these two facts together, we can say that
$$M_G^q(({\bf x}+{\bf y})([V]),U),\quad
M_G^q(({\bf x}\cdot{\bf y})([V]),U)$$
are also polynomials in
variables $\ld^n({\bf x}([V])),$ $\ld^n({\bf y}([V])),\, n\ge 1$, with integral coefficients.
Consequently, we can conclude that Eq. \eqref{poly identity} is an identity between polynomials in
$\ld^n({\bf x}([V])),\,\ld^n({\bf y}([V])),\,\, n\ge 1, [V] \preceq [U]$.
Obviously this identity makes sense for arbitrary special $\ld$-rings.
In the same way as above we can show that
\begin{equation*}
\begin{aligned}
&\tau^{q}({\bf x}+{\bf y})
=\tau^{q}({\bf x})+\tau^{q}({\bf y})
\end{aligned}
\end{equation*}
holds for arbitrary special $\ld$-rings. So we are done.
\qed

\begin{cor} \label{q-multiplicative of tau}
Let $A$ be a special $\ld$-ring equipped with $\mathbb Q$-algebra structure.
Then, for every $x,y \in A$, we have
$$\tau^G_q(qxy)=q\,\left(\tau^G_q(x)\cdot \tau^G_q(y)\right).$$
\end{cor}
\noindent{\bf Proof.}
Denote by $\varepsilon_{G}$ the vector determined by the condition
$\varepsilon_{G}([W])=\delta_{[G],[W]}.$
For every $x,y \in A$,
\begin{align*}
\tilde\varphi^q\left(q\,\left(\tau^G_q(x)\cdot \tau^G_q(y)\right)\right)
&=q \,\Phi^{q}(x\,\delta_G)\,\Phi^{q}(y\,\delta_G)\\
&=\dfrac 1q \,\Phi(q^2\,xy \,\delta_G)\\
&=\Phi^{q}(\,q\,xy \,\delta_G)\\
&=\tilde\varphi^{q}\left(\tau^G_q(qxy)\right)\,.
\end{align*}
The desired result follows from the injectiveness of $\varphi^q$.
\qed
\vskip 3mm
The next corollary is almost straightforward.
\begin{cor}
Let $A$ be a special $\ld$-ring equipped with $\mathbb Q$-algebra structure.
Then, as a ring,
$$\mathbb W^q_G(A) \stackrel{\text{\rm isomorphic} }{\cong}
Nr_G^q(A) \stackrel{\text{\rm isomorphic} }{\cong} {\widehat {Nr}}_G^q(A).$$
\end{cor}

\vskip 3mm
Finally, we will define natural transformations ${ v}_U^q,\, {f}_U^q$
on $\mathbb W^q_G$ by the transport of $q$-inductions and $q$-restrictions
via the map $\tau^q$.
First, in case $A=\mathbb Z$ we define
\begin{align*}
&v_U^q:=(\tau^q)^{-1}\circ \text{Ind}_U^{G}\circ \tau^q,\\
&f_U^q:=(\tau^q)^{-1}\circ q\text{-Res}_U^{G}\circ \tau^q.
\end{align*}

\begin{thm}\label{res and ind of q-witt burnside ring}
 Let $A=\mathbb Z$.
Then for every $[U]\in \mathcal O(G)$,
$v_U^q$ and $f_U^q$ are well-defined.
\end{thm}

\noindent {\bf Proof.}
For $\alpha \in \mathbb W_G(\mathbb Z)$ one has
\begin{align*}
&q\text{-Res}_U^G(\tau^q(\alpha))\\
&=\underset{[V]\in \mathcal O(G)}{{\sum}}
q\text{-Res}_U^G\cdot \text{Ind}_V^G(\tau^V_q(\alpha([V])))\\
&=\underset{[V] \in \mathcal O(G)}{{\sum}}\sum_{UgV\subseteq G}
\text{Ind}_{U\,\cap \, gVg^{-1}}^U \cdot
q\text{-Res}_{U \,\cap \, gVG^{-1}}^V(g)(\tau^V_q(\alpha([V])))\\
&=\underset{[V] \in \mathcal O(G)}{{\sum}}\sum_{UgV\subseteq G}
\text{Ind}_{U\,\cap \, gVg^{-1}}^U
(\tau^{U \, \cap \,gVg^{-1}}_q(\alpha([V])^{(V:U \, \cap \,gVg^{-1})})).
\end{align*}
Hence, it follows from Lemma \ref{Lem2} that
for any open subgroup $W$ of $U$
\begin{equation}\label{welldef of q witt-Burnside: restriction}
\begin{aligned}
&(\tau^q)^{-1}\circ q\text{-Res}_U^G\circ \tau^q(\alpha)([W])\\
&=\xi^{G,q}_{(W;W_1,\cdots,W_k;1,\cdots,1)}
(\alpha([V_1])^{(V_1;W_1)},\cdots,\alpha([V_k])^{(V_k,W_k)}),
\end{aligned}
\end{equation}
where $\xi^{G,q}_{(W;W_1,\cdots,W_k;1,\cdots,1)}$ is a polynomial with integral coefficients
(For complete information refer to Lemma \ref{Lem2}).
Similarly, one can show
$(\tau^q)^{-1}\circ \text{Ind}_U^G\circ \tau^q(\alpha)([W])$
is a polynomial with integral coefficients in those $\alpha([V])$'s,
where $V$ ranges over the set of open subgroups of $U$
to which $W$ is sub-conjugate in $G$.
\qed
\vskip 3mm
For an arbitrary commutative ring $A$,
let us define restriction $f_U^q: \mathbb W^q_G(A) \to \mathbb W^q_U(A)$
using the polynomials appearing in Eq. \eqref{welldef of q witt-Burnside: restriction}.
This will give us a natural transformation $f_U^q:\mathbb W^q_G \to \mathbb W^q_U$ satisfying
$$\tau^q\circ f_U^q= q\text{-Res}_U^{G}\circ \tau^q.$$
Similarly, one can define a natural transformation
$v_U^q: \mathbb W^q_U \to \mathbb W^q_G$ such that
$$\tau^q\circ v_U^q= \text{Ind}_U^{G}\circ \tau^q.$$

\section{Lenart's conjecture and proof of lemmas and theorems}\label{Proofs}
\subsection{Lenart's conjecture}
In \cite[page 731]{L}) Lenart proposed a conjecture on $q$-restriction
defined on $Nr_G^q(A)$ in case where $G=\hat C$ (see Section \ref{Necklace}).
It can be stated as follows.
\vskip 3mm
\noindent{\bf Conjecture}:\\
{\it Set
$$f^q_r:=\mathcal F_{\hat
C^r}^{\hat C}$$ and
\begin{align*}
M^q(x,n)=\sum_{d|n}\mu^q(n,d)\,q^{d-1}\,x^d.
\end{align*}
Then the $n$-th component of $f^q_r\,\left(M^q(x,n)_{n\ge 1}\right)$ is given by
\begin{equation*}
\sum_{d|n}Q_{r,n,d}(q)M^q(x^r,d) \quad (\text { in } \mathbb Q[x,q]),
\end{equation*}
where $Q_{r,n,d}(q) \in \mathbb Q[q]$ are numerical polynomials.}
\vskip 3mm
\noindent{\bf Proof.}
Eq. \eqref{relation of q-restriction} implies that
\begin{equation*}
f^q_rM^q(x)=M^q(q^{r-1}x^r).
\end{equation*}
If $r=1$, there is nothing to prove since $f^q_1$ is the identity map.
For $r\ge 2$
\begin{align*}
M^q(q^{r-1}x^r)
&=M^q(q \cdot q^{r-2}\cdot x^r)\\
&=\,q (M^q(q^{r-2})\cdot M^q(x^r)),
\end{align*}
where the last equality follows from
Corollary \ref{q-multiplicative of tau} or \cite [Proposition 5.15]{L}.
Since
\begin{align*}
&q\left(M^q(q^{r-2})\cdot  M^q(x^r)\right)(n)\\
&=q\sum_{[i,j]\,|\,n} (i,j)P^n_{i,j}(q)M^q(q^{r-2},i)M^q(x^r,j)\\
&=q\sum_{j|n}
\left(\sum_{[i,j]\,|\,n} (i,j)P^n_{i,j}(q)M^q(q^{r-2},i)\right)M^q(x^r,j),
\end{align*}
we obtain
\begin{equation*}
Q_{r,n,d}(q)=q\left(\sum_{[i,d]\,|\,n} (i,d)P^n_{i,d}(q)M^q(q^{r-2},i)\right).
\end{equation*}
But, it was shown in \cite {L} that $P^n_{i,d}(q)$ is a numerical polynomial.
And $M^q(q^{r-2},i)$ is also a numerical polynomial in view of
Theorem \ref{expression of q-neck associ profinite group}.
Thus, we can conclude that $Q_{r,n,d}(q)$ are numerical polynomials.
\qed
\subsection{Proof of Theorem \ref{well of q-restric}}\label{proof of q-rest}
Recall that the matrix representing $q\text{-}{\rm Res}_U^G$
is given by
\begin{align*}
d_{[V],[W]}(q)=
\begin{cases}
Q_{V,W}(q)\Psi^{(W:V)}&\text{ if } [W]\preceq [V] \text{ in }
\mathcal O(G),\\
0&\text{ otherwise,}
\end{cases}
\end{align*}
where
\begin{equation*}
Q_{V,W}(q)=\sum_{[S]\in \mathcal O(U)
\atop {[W]\preceq [S]\text{ in }\mathcal O(G)
\atop  [S]\preceq
[V]\text{ in }\mathcal O(U)}}
\mu_U^q([V],[S])\varphi_S(G/W)q^{(W:S)-1}
\end{equation*}
Refer to Section \ref{on indiction and restriction}.
It is straightforward that $Q_{V,W}(q)$ is given by
\begin{align*}
q\text{\rm -Res}_U^G (\varepsilon_{[W]})([V]),
\end{align*}
where $\varepsilon_{[W]}$ represents the vector determined by the condition
\begin{equation}\label{standard basis}
\varepsilon_{[W]}([W'])=\delta_{[W],[W']}.
\end{equation}
By Theorem \ref{ringhomo of q-teichmuller*}
we can write $\varepsilon_{[W]}$ as
$${\underset{[U]\in \mathcal O(G)}{\sum}}\,\text{Ind}_U^{G}
\circ \tau^V_q(\alpha([V]))$$
for some $\alpha \in A$.
Now, from $$\mathbb W_G^q(\mathbb Z)\cong Nr_G^q(\mathbb Z),$$
it follows that
$\alpha([V])\in \mathbb Z$ for all $[V]\in \mathcal O(G)$.
Then, by Proposition \ref{Mackey} (b),
\begin{align*}
&q\text{\rm -Res}_U^G \circ
\left({\underset{[V]\in \mathcal O(G)}{\sum}}\,\text{Ind}_V^{G}
\circ \tau^V_q(\alpha([V]))\right)\\
&=\sum_{UgV\subseteq G}{\underset{[V]\in \mathcal O(G)}{\sum}}
\text{\rm Ind}_{U\cap gVg^{-1}}^U
\circ
q\text{\rm -Res}_{U\cap gVg^{-1}}^V (g)
\left(\tau^V_q(\alpha([V]))\right)\\
&=\sum_{UgV\subseteq G}{\underset{[V]\in \mathcal O(G)}{\sum}}
\text{\rm Ind}_{U\cap gVg^{-1}}^U
(\tau_q^{U\cap  gVg^{-1}}(q^{(V:U\cap  gVg^{-1})-1}\alpha([V])^{(V:U\cap  gVg^{-1})}).
\end{align*}
The last equality follows from Eq. \eqref{relation of q-restriction}.
Finally, our assertion follows from the fact that
$M^q_G(x,V)$ is a numerical polynomial
for every open subgroup $V$ of $G$.
\qed
\subsection{Proof of Lemma \ref{form of product of necklace ring**}}
\begin{lem}\label{lem for numericality of necklace}
Let  $n\in \mathbb Z_{>0}$ and $q\in \mathbb Z$.
Write
\begin{align*}
n=p_1^{a_1}\cdots p_r^{a_r}n',\quad
q=p_1^{b_1}\cdots p_r^{b_r}q',
\end{align*}
where $p_i$'s are primes, $a_i, b_i >0$, and $ (n',p_i)=(q',p_i)=(n',q')=1$
for all $i$.
Then,
$p_i^{a_i+b_i}$ divides $\frac nd q^{d}$
for all $1\le i \le r$ and for all $d$ dividing $n$.
\end{lem}
\noindent{\bf Proof.}
Write $d=p_1^{c_1}\cdots p_r^{c_r}d'$ where $(p_i,d')=1$ for all $1 \le i \le r.$
Then $p_i^{a_i-c_i}p_i^{b_ip_i^{c_i}}$ divides $\frac nd q^{d}.$
But since $a_i-c_i+b_ip_i^{c_i} \ge a_i+b_i$ for all $1\le i \le r,$
we get the desired result.
\qed

\begin{lem}\label{numericality of necklace}
Let $q$ be a non-zero integer.
For an open subgroup $V$ of $G$ let
\begin{equation*}
M_G^q(x,V):=\sum_{[W]\preceq [V]}\mu^q_{G}([V],[W])\,q^{(G:W)-1}\,x^{(G:W)}\,.
\end{equation*}
If $x=q^m$ for some positive integer $m$,
then $\dfrac 1q \,M_G^q(x,V)$ is a numerical polynomial in $q$.
\end{lem}

\noindent{\bf Proof.}
In view of Lemma \ref{relation of hom and necklace*}
and Theorem \ref{expression of q-neck associ profinite group}
we have
\begin{equation*}
q^{(G:V)-1}\,q^{m(G:V)}
=\sum_{[W]\preceq [V]}\varphi_V(G/W)q^{(W:V)-1}M^q_G(q^m,W),
\end{equation*}
equivalently
\begin{equation}\label{diveded by q?}
\varphi_V(G/V)M^q_G(q^m,V)=q^{(G:V)-1}\,q^{m(G:V)}
-\sum_{[W]\preceq [V] \atop [W]\ne [V]}\varphi_V(G/W)\,q^{(W:V)}\,\frac 1q M^q_G(q^m,W).
\end{equation}
Note that $M^q_G(q^m,G)=q^m$.
Hence, we can use an induction on the index $(G:V)$.
Assume that our assertion holds for all $W$'s satisfying $(G:W)>(G:V)$.
Write
\begin{align*}
&\varphi_V(G/V)=p_1^{a_1}\cdots p_r^{a_r}\alpha'\\
&q=p_1^{b_1}\cdots p_r^{b_r}q',
\end{align*}
where $p_i$'s are primes, $a_i, b_i >0$, and $ (\alpha',q')=1$.
Clearly $q'\,|\,M^q_G(q^m,V)$ since $q$ divides the right hand side of
Eq. \eqref{diveded by q?}.
Recall that
Corollary \ref{two same formulas}
implies that
$$\varphi_V(G/W)
=\sum_{[V_i]\in \mathcal O(W)\atop [V_i]
=[V]\,\,\, \text{\rm in }\mathcal O(G)}[N_G(V_i):N_W(V_i)].$$
Since $[N_G(V_i):N_W(V_i)]\,|\, \varphi_V(G/V)$ and
$(W:V) \ge [N_W(V_i): V_i]$
we can obtain the desired result by applying Lemma \ref{lem for numericality of necklace}.
\qed

\begin{lem}
For every $[V],[W] \in \mathcal O(G)$
and for every $\alpha, \beta \in \mathbb Z$
one has the $q$-modified Mackey formula
\begin{equation}\label{q-Mackey formula}
\begin{aligned}
&q (\text{\rm Ind}_V^G(\tau^V_q(\alpha))\cdot \text{\rm Ind}_W^G(\tau^W_q(\beta)))\\
&=\sum_{VgW \subseteq G}\text{\rm Ind}_{V\cap gWg^{-1}}^G \circ \tau^{V\cap gWg^{-1}}_q\\
&\quad \quad \left (q^{-1}(q^{(V:V\cap gWg^{-1})} \alpha^{(V:V\cap gWg^{-1})})
(q^{(W:V\cap gWg^{-1})} \beta^{(W:V\cap gWg^{-1})})\right ).
\end{aligned}
\end{equation}
\end{lem}

\noindent{\bf Proof.}
If $q=0$, then it is trivial.
So we assume that $q\ne 0$.
We show that $\tilde \varphi^q_U$, applied to both sides
of the above equation, yields the same number for all open subgroups
$U$ in $G$.
Observe
\begin{align*}
&q\left( \tilde \varphi^q_U
\left(\text{Ind}_V^G(\tau^V_q(\alpha))
\cdot \text{Ind}_W^G(\tau^W_q(\beta)\right)\right)\\
&=q \left(\tilde\varphi^q_U
\left(\text{Ind}_V^G(\tau^V_q(\alpha)\right)
\cdot \tilde \varphi^q_U\left(\text{Ind}_W^G
(\tau^W_q(\beta)\right)\right)\\
&=q \left(\varphi_U(G/V)\,q^{(V:U)-1}\alpha^{(V:U)}
\cdot \varphi_U(G/W)\,q^{(W:U)-1}\beta^{(W:U)}\right)
\quad \text{ (by Theorem \ref{ringhomo of q-teichmuller} (b)})\\
&=q\left(\varphi_U(G/V \times G/W)\, (q^{(V:U)-1}\alpha^{(V:U)})
\,(q^{(W:U)-1}\beta^{(W:U)}) \right)\\
&=\sum_{VgW\subseteq G}\varphi_U(G/V\cap gWg^{-1})\, q^{(V\cap gWg^{-1}:U)-1}\\
&\quad \times
\left(q^{-1}(q^{(V:V\cap gWg^{-1})}\alpha^{(V:V\cap gWg^{-1})})
(q^{(W:V\cap gWg^{-1})-1}\beta^{(W:V\cap gWg^{-1})})\right)^{(V\cap gWg^{-1}:U)}\\
&=\tilde\varphi^q_U \left(\sum_{VgW\subseteq G}
\text{Ind}_{V\cap gWg^{-1}}^G \circ \tau^{(V\cap gWg^{-1}:U)}_q
\left(q^{-1}(q^{(V:V\cap gWg^{-1})}\alpha^{(V:V\cap gWg^{-1})})\right.\right.\\
& \qquad  \left. \left.
(q^{(W:V\cap gWg^{-1})}\beta^{(W:V\cap gWg^{-1})})
\right)\right).
\end{align*}
\qed
\vskip 3mm

\noindent{\bf Main Proof.}
In case where $q=0$, our statement will be trivial since $\tilde\varphi^0$ is the identity map.
Therefore, we assume that $q$ is not zero. In view of
Lemma \ref{form of product of necklace ring} we obtain that for
every $[V]\in \mathcal O(G)$
\begin{equation*}
\varepsilon_{[V]}\cdot \varepsilon_{[W]}([U])=
\begin{cases}
P_{V,W}^U(q)& \text{ if } [V],[W]\preceq [U]\\
0 &\text{ otherwise.}
\end{cases}
\end{equation*}
(see Eq. \eqref{standard basis}).
Note that
$$\tau^q(\varepsilon_{[V]})=\varepsilon_{[V]}\quad
\text{ and }\quad
\tau^q(\varepsilon_{[W]})=\varepsilon_{[W]},
$$
which can be shown by comparing the image of each side for $\tilde\varphi$.
Eq. \eqref{q-Mackey formula} implies that if $[V],[W]\preceq [U]$, then
\begin{equation}\label{concrete form of coefficient}
\begin{aligned}
&P_{V,W}^U(q)\\
&=\dfrac 1q\sum_{VgW \subseteq G}
\text{\rm Ind}_{V\cap gWg^{-1}}^G \tau_q^{V\cap gWg^{-1}}
\left(q^{(V:V\cap gWg^{-1})+(W:V\cap gWg^{-1})-1}\right)([U]).
\end{aligned}
\end{equation}
From Lemma \ref{numericality of necklace}
it follows that
$P_{V,W}^U(q)$ is a numerical polynomial in $q$.
\qed

\begin{exm}\label{explicit of p}{\rm
If $G$ is abelian, then
\begin{equation*}
\begin{aligned}
P_{V,W}^U(q)
=\dfrac 1q(G:V+W)
M^q_{V\cap W}
\left(q^{(V:V\cap W)+(W:V\cap W)-1}, U\right).
\end{aligned}
\end{equation*}
If $G=\hat C$ and $i,j\,|\,n$, then
$$P_{\hat C^i,\hat C^j}^{\hat C^n}(q)=\frac {(i,j)}{q} M^q\left(q^{\frac {i+j}{(i,j)}-1}, \frac {n}{[i,j]}\right).$$
}\end{exm}
\begin{rem}{\rm
For all $[U],[V],[W]$'s satisfying $[V],[W] \preceq [U]$,
$P_{V,W}^U(q)$ can also be computed recursively via
the following formula:
\begin{equation*}
\sum_{[V],[W]\preceq [Z]\preceq [U]}
\varphi_U(G/Z)q^{(Z:U)-1}P_{V,W}^Z(q)=\varphi_U(G/V)\varphi_U(G/W)
\,q^{(V:U)+(W:U)-2}.
\end{equation*}
Indeed, this identity follows
from the fact that $\tilde \varphi^q$ is a ring homomorphism.
In particular, if $U=Z(g,V,W)$, then
\begin{align*}
P_{V,W}^U(q)=\frac {1}{(N_G(U):U)}\varphi_U(G/V)\varphi_U(G/W)
\,q^{(V:U)+(W:U)-2}.
\end{align*}
}\end{rem}
\subsection{Proof of Theorem \ref{q-deformation of witt-Burn}}
\begin{lem}{\rm (cf. \cite [Lemma (3.2.2)]{DS2})}
\label{Lem 8.1}
With the notation in \cite [Lemma (3.2.2)]{DS2},
we obtain that for any $\alpha,\beta \in R$
\begin{equation}\label{expansion of addition}
\tau^G_q(\alpha+\beta)
=\sum_{G\cdot A \in G\setminus \mathfrak U(G)}
{\rm Ind}_{U_A}^{G}\left(\tau^{U_A}_q
\left(q^{-1}(q^{i_A }\alpha^{i_A}\cdot q^{i_{G-A}}\beta^{i_{G-A}})\right)\right).
\end{equation}
\end{lem}
\noindent{\bf Proof.}
First we assume that $R$ is torsion-free.
For every open subgroup $U$ of $G$,
if we take $\tilde \varphi_U^q$ on the right side of
Eq.\eqref{expansion of addition}, one has
\begin{align*}
&\tilde \varphi_U^q
\left(\sum_{G\cdot A \in G\setminus \mathfrak U(G)}
\text{\rm Ind}_{U_A}^{G}\left(\tau^{U_A}_q
\left(q^{-1} (q^{i_{A}}\alpha^{i_A}\cdot q^{i_{G-A}}\beta^{i_{G-A}})\right)\right)\right)\\
&=\sum_{G\cdot A \in G\setminus \mathfrak U(G)\atop U\lesssim U_A}
\tilde \varphi_U(G/U_A)\,q^{(U_A:U)-1}
\left(q^{-1}( q^{i_A }\alpha^{i_A}
\cdot q^{i_{G-A}}\beta^{i_{G-A}})\right)^{(U_A:U)}
\text{ (by  Theorem \ref{ringhomo of q-teichmuller} (b)
)}\\
&=q^{\sharp (G/U)-1} \sum_{A \in \mathfrak U(G),\, U \lesssim U_A}
\alpha^{\sharp (A/U)}\cdot \beta^{\sharp (G-A)/U}\\
&=q^{\sharp (G/U)-1}(\alpha+\beta)^{\sharp (G/U)}\\
&=\tilde \varphi_U^q(\tau^G_q(\alpha+\beta)).
\end{align*}
Since $\tilde \varphi_U^q$ is injective, we have the desired result.
In case $R$ is not torsion-free,
$\tilde \varphi_U^q$ is no longer injective.
However we note that the $U$-th component appearing in
Eq. \eqref{expansion of addition} is an integral polynomial in
$\ld^k(\alpha_V)$'s and $\ld^l(\beta_W)$'s for
$1 \le k,l \le  [G:U]$ and $U \lesssim V,W \leqslant G$
for an arbitrary open subgroup $U$ of $G$.
This implies that Eq. \eqref{expansion of addition}
holds regardless of torsion.
\qed

\begin{lem}{\rm (cf. \cite[Lemma (3.2.5)]{DS2})} \label{Lem2}
For some $k\in \mathbb N$ let $V_1,\cdots,V_k \leqslant G$
be a sequence of open subgroups of $G$.
Then for every open subgroups $U \leqslant G$ and every sequence
$\varepsilon_1,\cdots,\varepsilon_k \in \{\pm 1\}$
there exists a unique polynomial
$\xi_U^q
=\xi^{G,q}_{(U;V_1,\cdots,V_k;\varepsilon_1,\cdots,\varepsilon_k)}
=\xi_U^q(x_1,\cdots,x_k)\in \mathbb Z[x_1,\cdots,x_k]$
such that for all $\alpha_1,\cdots,\alpha_k \in \mathbb Z$
one has
$$(\tau^q)^{-1}(\sum_{i=1}^k
\varepsilon_i\cdot \text{\rm Ind}_{V_i}^G(\tau^{V_i}_q(\alpha_i))(U)
=\xi_U^q(\alpha_1,\cdots,\alpha_k).$$
\end{lem}

\noindent{\bf Proof.}
The proof can be done in the exactly same way of that of
Lemma (3.2.5) \cite{DS2}.
Without loss of generality one may assume
$$U\leqslant \bigcap_{i=1}^k V_i.$$
If $\varepsilon_1=\varepsilon_2=\cdots=\varepsilon_k=1$ and if $V_i$
is not conjugate to $V_j$ for $i\ne j$, then
$\sum_{i=1}^k \varepsilon_i
\text{Ind}_{V_i}^G(\tau^{V_i}_q(\alpha_i))=\tau^q(\alpha)$
for $\alpha \in \mathbb W_G^q(\mathbb Z)$ with
$[V]=[V_j]$ and $\alpha([V])=0$ if $V$ is not conjugate to either of the
$V_1,\cdots, V_k$.
So in this case we are done :
$\xi_U^q(\alpha_1,\cdots,\alpha_k)$ equals $\alpha_{[U]}$.
Hence, we may use triple induction, first with respect to
\begin{align*}
&m_1=m_1(U;v_1,\cdots,V_k;\varepsilon_1,\cdots, \varepsilon_k):
=\text{max}((V_i:U)\,|\,\varepsilon_i=-1 \\
&\qquad\text{ or there exists some }j\ne i \text{ with } V_j\sim V_i),
\end{align*}
then with respect to
$$m_2:=\sharp \, \{ \, i\, | \,(V_i:U)=m_1 \text{ and } \varepsilon_i=-1\},$$
and then with respect to
\begin{align*}
m_3:=\sharp \, \{\, i\, | \,(V_i:U)=m_1 \text{ and there exists some }
j\ne i \text{ with }V_j \sim V_i\}.
\end{align*}
In case where $m_1=0$, the induction hypothesis holds in view of the above
remark.
In case where $m_1>0$, we have either $m_2>0$ or $m_3>0$.
In case $m_2>0$, say $(V_1:U)=m_1$ and $\varepsilon_1=-1$,
we may use Eq. \eqref{Lem 8.1}
with $G=V_1, \alpha=-\alpha_1$, $\beta=+\alpha$ to conclude that
\begin{equation}
0=\sum_{V_1\cdot A\in V_1 \setminus \mathfrak U(V_1)}
{\rm Ind}_{U_A}^{V_1}(\tau^{U_A}_q
((-1)^{i_A}q^{(V_1:U_A)-1}\alpha_1^{(V_1:U_A)})).
\end{equation}
Therefore, considering the two special summands
$A=\emptyset$ and $A=V_1$ and putting
$\mathfrak U_0(V_1):=\{A\in \mathfrak U(V_1):
A\ne \emptyset \text { and } A\ne V_1\}$, one gets
\begin{equation*}
-\tau^{U_A}_q(\alpha_1)
=\tau^{U_A}_q(-\alpha_1)
+\sum_{V_1\cdot A\in V_1 \setminus \mathfrak U_0(V_1)}
{\rm Ind}_{U_A}^{V_1}(\tau^{U_A}_q
((-1)^{i_A}\cdot q^{(V_1:U_A)-1}\alpha_1^{(V_1:U_A)})).
\end{equation*}
Hence, if $A_{k+1},A_{k+2},\cdots, A_{k'}\in \mathfrak U_0(V_1)$
denote representatives of the $V_1$-orbits
$V_1\cdot A \subseteq \mathfrak U_0(V_1)$ with
$$U\leqslant V_{k+1}:=U_{A_{k+1}},\, V_{k+2}
:=U_{A_{k+2}},\cdots,\, V_{k'}:=U_{A_{k'}}\lneqq V_1$$
and
if we put $\varepsilon_{k+1}=\cdots=\varepsilon_{k'}:=1$
and
$$\alpha_{k+1}:=(-1)^{i_{A_{k+1}}}q^{(V_1:V_{k+1})-1}\cdot \alpha_1^{(V_1:V_{k+1})},
\cdots, \alpha_{k'}:=(-1)^{i_{A_{k'}}}q^{(V_1:V_{k'})-1}\cdot \alpha_1^{(V_1:V_{k'})},$$
then we have
\begin{align*}
&{\xi_U^q}_{(U;V_1,\cdots,V_k;-1,\varepsilon_2,\cdots,\varepsilon_k)}
(\alpha_1,\alpha_2,\cdots,\alpha_k)\\
&={\xi_U^q}_{(U;V_1,\cdots,V_{k'};1,\varepsilon_2,\cdots,\varepsilon_{k'})}
(-\alpha_1,\alpha_2,\cdots,\alpha_{k'}),
\end{align*}
so the result follows by induction.

Similarly, if $m_2=0,$ but $m_3>0$, say $V_1=V_2$,
then we may use Eq. \eqref{Lem 8.1} once more with
$G=V_1, \alpha=\alpha_1,$ and $\beta=\alpha_2$
to conclude that
\begin{align*}
&\tau^{U_A}_q(\alpha_1+\alpha_2)\\
&=\sum_{V_1\cdot A\in V_1 \setminus \mathfrak U(V_1)}
{\rm Ind}_{U_A}^{V_1}(\tau^{U_A}_q
(q^{(V_1:U_A)-1}\alpha^{i_A}\alpha_2^{i_{V_1-A}}))\\
&=\tau^{U_A}_q(\alpha_1)+\tau^{U_A}_q(\alpha_2)
+\sum_{V_1\cdot A\in V_1 \setminus \mathfrak U_0(V_1)}
{\rm Ind}_{U_A}^{V_1}(\tau^{U_A}_q
(q^{(V_1:U_A)-1}\alpha_1^{i_A}\alpha_1^{i_{V_1-A}})).
\end{align*}
So with $V_{k+1},\cdots,V_{k'}$ as above, but
$\varepsilon_{k+1}=\cdots=\varepsilon_{k'}=-1$ and
$$\alpha_{k+1}:=q^{(V_1:U_A)-1}\alpha_1^{i_{A_{k+1}}}\alpha_2^{i_{V_1-A_{k+1}}},\,
\cdots,\, \alpha_{k'}:=q^{(V_1:U_A)-1}\alpha_1^{i_{A_{k'}}}\alpha_2^{i_{V_1-A_{k'}}},$$
we get
\begin{align*}
&{\xi_U^q}_{(U;V_1,\cdots,V_k;\varepsilon_1,\varepsilon_2,\cdots,\varepsilon_k)}
(\alpha_1,\alpha_2,\cdots,\alpha_k)\\
&={\xi_U^q}_{(U;V_2,\cdots,V_{k'};\varepsilon_2,\cdots,\varepsilon_{k'})}
(\alpha_1+\alpha_2,\alpha_3,\cdots,\alpha_{k'}).
\end{align*}
So again our result follows by induction
since $A\in \mathfrak U_0(V_1)$ implies that
$U_A\precneqq V_1$.
\qed
\vskip 3mm

\noindent{\bf Main Proof.}
First, let us consider the equation
\begin{equation*}
\Phi^q({\bf x})=\Phi^q(\alpha)+\Phi^q(\beta).
\end{equation*}
Applying the identity \eqref{relatio of phi and q-phi}
yields that
\begin{align*}
\Phi^q(s^q)
&=\frac 1q \Phi(q\alpha)+\frac 1q \Phi(q\beta)
\quad (\text{ in } {\rm Gh}(G,\mathbb Z)) \\
&=\frac 1q \,\Phi(q\alpha \,+\,q\beta)
\quad (\text{ in } {\mathbb W}_G(\mathbb Z))\\
&=\Phi^q \left(\frac 1q(q\alpha \,+\,q\beta)\right)
\quad (\text{ in }{\mathbb W}_G(\mathbb Z)).
\end{align*}
Thus, we get
$$s_U^q=\frac 1q \,s_U(\,qx_V,qy_V\,:\, [V]\preceq [U]\,).$$
Clearly $s_U^q$ is a polynomial in
$\mathbb Z\,[\,{\bf x}([V]),{\bf y}([V])\,:\,[V]\preceq [U]\,]$
since $s_U$ has no constant term.
Similarly,
\begin{align*}
-\Phi^q({\bf a})
&=-\frac 1q \Phi(q{\bf a})\\
&=\frac 1q \, \Phi(-q{\bf a}) \quad (`-' \text{means the inverse of $+$ in }\mathbb W(\mathbb Z))\\
&=\Phi^q\left( \frac 1q (-q{\bf a})\right).
\end{align*}
From this it follows that
$$\iota_{[U]}^q=\frac 1q \iota_{[U]}\left((-q{\bf a})([V]) \,:\, [V] \preceq [U]\,\right).$$
Clearly it has integer coefficients
since $\iota_{[U]}$ has integer coefficients and no constant term.

To compute $p^q_U$ we use Eq. \eqref {q-Mackey formula}.
First we choose a system $s_1,s_2,\cdots,s_h$ of representatives
of the $G$-orbits in
$$S:=\underset {1 \le i,j \le k }{\dot\bigcup} G/V_i \times G/V_j .$$
Next, we put $W_r:=G_{s_r}$ and
$$p_r^q(x_{V_1},y_{V_1},\cdots,x_{V_k},y_{V_k})
:=q^{(V_i:W_r)-1}x_i^{(V_i:W_r)}q^{(V_j:W_r)-1}y_i^{(V_j:W_r)}$$
in case $s_r=(g_rV_i,g'_{r}V_j) \in G/V_i \times G/V_j \subseteq S.$
Using these conventions and Eq. \eqref{q-Mackey formula}, we get the equation
\begin{align*}
&p_U^{q}(x_{V_1},y_{V_1},\cdots,x_{V_k},y_{V_k})\\
&=\xi^{G,q}_{(U;W_1,\cdots,W_h;1;\cdots;1)}
(p_1^q(x_{V_1}\cdots,y_{V_k}),\cdots,p_h^q(x_{V_1}\cdots,y_{V_k})).
\end{align*}
So, $p_U^{q}$ is also in
$\mathbb Z\,[\,{\bf x}([V]),{\bf y}([V])\,:\,[V]\preceq [U]\,].$
\qed
\section{Classification $\mathbb W_G^q$ up to strict natural isomorphism}
In this section, we assume that $G$ is an abelian profinite group.
The aim of this section is to classify $\mathbb W_G^q$
up to strict natural isomorphism as $q$ ranges over the set of
integers.
To begin with, we introduce prerequisites.
Given an integer $q$, we denote by $D(q)$
the set of divisors of $q$, and by $D^{{\rm pr}}(q)$ the set of
prime divisors of $q$, respectively.
Conventionally, $D(0)$ will
denote the set of positive integers $\mathbb N$, and $D^{{\rm
pr}}(0)$ the set of all primes in $\mathbb N$.

\begin{df}
Let $q$ and $r$ be arbitrary integers.

{\rm (a)} Given a commutative ring $A$,
$\mathbb W_G^q(A)$ is said to be {\it strictly-isomorphic} to $\mathbb W_G^r(A)$
if there exists a ring isomorphism,
say $\tau_q^r: \mathbb W_G^q(A) \to \mathbb W_G^r(A)$, satisfying
$\Phi_G^q=\Phi_G^r\circ \tau_q^r$.
In this case, $\tau_q^r$ is called a {\it strict-isomorphism}.

{\rm(b)}
$\mathbb W_G^q$ is said to be {\it strictly-isomorphic} to $\mathbb W_G^r$
if there exists a natural isomorphism,
say $\tau_q^r: \mathbb W_G^q\to \mathbb W_G^r$, satisfying
$\Phi_G^q=\Phi_G^r\circ \tau_q^r$.
In this case, $\tau_q^r$ is called a {\it strict natural isomorphism}.
\end{df}

Denote by ${\rm D}^{\rm pr}(G)$ the set of prime divisors
of each of $\{(G:U)
:U\in \mathcal O(G)\}.$
And we let
\begin{align*}
&D^{\rm{pr}}(q)\cap D^{\rm pr}(G)=\{p_1,\cdots,p_k ; c_1,\cdots,c_s \},\\
&D^{\rm{pr}}(r)\cap D^{\rm pr}(G)=\{p_1,\cdots,p_k ; d_1,\cdots,d_t \}.
\end{align*}
That is, $p_i$'s are primes in
$D^{\rm{pr}}(q)\cap D^{\rm{pr}}(r)\cap D^{\rm pr}(G).$
\begin{thm}\label{main theorem:Witt- Burnside}
Let $A$ be a commutative ring with identity.
And, let $q,r$ be arbitrary integers and $G$ an abelian profinite group.
Then,
there exists a unique strict-isomorphism between
$\mathbb W_G^q(A)$ and $\mathbb W_G^r(A)$
if and only if $A$ is a
$\mathbb Z\,[\,\frac {1}{c_i},
\frac {1}{d_j} \,:\, 1\le i \le s, 1\le j \le t\,]$-algebra.
\end{thm}
The above result will be proved in the following steps.
First, let
$$R=\mathbb Q[X_U\,:\, U\in \mathcal O(G)],$$
and then consider the set of equations
arising from
$$\Phi_G^q({\bf X})=\Phi_G^r({\bf Y}).$$
More precisely, this set consists of the following identities
\begin{equation}\label{r-phi=q-phi:Witt-Burnside}
\begin{aligned}
\sum_{V\in \mathcal O(G) \atop
U\subseteq V\subseteq G}
(G:V)q^{(V:U)-1}
X_{V}^{(V:U)}
=\sum_{V\in \mathcal O(G) \atop
U\subseteq V\subseteq G}(G:V)r^{(V:U)-1}
Y_{V}^{(V:U)}
\end{aligned}
\end{equation}
for all $U\in \mathcal O(G)$.
Here, ${\bf X}$ and ${\bf Y}$ represents
the vectors $(X_{U})_{U\in \mathcal O(G)}$ and
$(Y_{U})_{U\in \mathcal O(G)}$, respectively.
Denote by $\mathbb Z_{G}$ the commutative ring
$\mathbb Z\left[\frac 1p: p \in D^{\rm pr}(G)\right]$
(see Eq. \eqref{defn of local}).
\begin{lem}\label{key lemma: Witt-Burnside}
Let $G$ be an arbitrary profinite group.
Then, for every $U\in \mathcal O(G)$, it holds that
$$Y_{U}-X_{U}\in
\left(\mathbb Z_G \cap \mathbb Z\left[\frac 1q,\frac 1r\right]\right)
\left[X_{S}\,:\, U\subsetneqq S \subseteq G, \,\,\,{\rm and }\,\, S\in \mathcal O(G) \right].$$
\end{lem}

\noindent{\bf Proof.}
Given $U\in \mathcal O(G)$ we would like to express $Y_{U}$
as a polynomial in $X_{V}$'s inductively.
In view of Eq. \eqref{r-phi=q-phi:Witt-Burnside},
$Y_{G}=X_{G}$.
Now, we assume that
$$Y_{V}-X_{V} \in
\mathbb Z_G \,[X_{S}\,:\, V\subsetneqq S \subseteq G, \,\,\,{\rm and }\,\, S\in \mathcal O(G)\,]$$
for all $V\in \mathcal O(G)$ with $U\subsetneqq V \subseteq G$.
Transform Eq. \eqref{r-phi=q-phi:Witt-Burnside}
into the form
\begin{equation*}
\begin{aligned}
Y_{U}-X_{U}=\sum_{V\in \mathcal O(G)\atop
U\subsetneqq V \subseteq G}
\frac {(q^{(V:U)-1} X_{V}^{(V:U)}-r^{(V:U)-1}Y_{V}^{(V:U)})}
{(V:U)}.
\end{aligned}
\end{equation*}
The induction hypothesis immediately implies that
\begin{equation*}\label{conclusion 1}
Y_{U}-X_{U}\in
\mathbb Z_G \,[X_{S}\,:\, U\subsetneqq S \subseteq G, \,\,\,{\rm and }\,\, S\in \mathcal O(G)\,].
\end{equation*}
Next, let us show that the coefficients are in $\mathbb Z\left[\frac 1q,\frac 1r\right].$
Let us first assume that $q$ and $r$ are nonzero.
Multiply
 $qr/(q,r)$ to both sides of Eq. \eqref{r-phi=q-phi:Witt-Burnside}
to obtain the identity
$$\Phi_G\left((qX_{U})_{U\in \mathcal O(G)}\right)
=\Phi_G\left(\frac {q}{(q,r)}\cdot (rY_{U})_{U\in \mathcal O(G)}\right).$$
It follows from the injectiveness of $\Phi_G$
that
\begin{equation*}\label{q,r: coef of Witt-Burnside}
\frac {r}{(q,r)}\cdot (qX_{U})_{U\in \mathcal O(G)}
=\frac {q}{(q,r)}\cdot (rY_{U})_{U\in \mathcal O(G)}.
\end{equation*}
It is not difficult to show that, for all $U\in [G]$,
the $U$-th coordinates of
$r/(q,r)\cdot (qY_{U})_{U\in \mathcal O(G)}$ and $q/(q,r)\cdot (qX_{U})_{U\in \mathcal O(G)}$
are of the form
\begin{equation}\label{second induction}
\frac {rq}{(q,r)}Y_{U}+ \text{ a polynomial contained in }
\mathbb Z \,[qY_{S}\,:\,U\subsetneqq S \subseteq G, \,\,\,{\rm and }\,\, S\in \mathcal O(G)\,]
\end{equation}
and
\begin{equation}\label{second induction*}
\frac {rq}{(q,r)}X_{U}+ \text{ a polynomial contained in }
\mathbb Z \,[qX_{S}\,:\, U\subsetneqq S \subseteq G, \,\,\,{\rm and }\,\, S\in \mathcal O(G)\,],
\end{equation}
respectively.
Assume that
$$Y_{V}-X_{V}\in \mathbb Z\left[\frac 1q,\frac 1r\right]
\,[X_{S}\,:\, V\subsetneqq S \subseteq G, \,\,\,{\rm and }\,\, S\in \mathcal O(G)\,]$$
for all $V \in \mathcal O(G)$ with $U \subsetneqq V \subseteq G.$
Applying the induction hypothesis above to
Eq. \eqref{second induction} and Eq. \eqref{second induction*} yields
\begin{equation*}\label{conclusion 2}
Y_{U}-X_{U} \in \mathbb Z\left[\frac 1q,\frac 1r\right]
\,[X_{S}\,:\, U\subsetneqq S \subseteq G, \,\,\,{\rm and }\,\, S\in \mathcal O(G) \,].
\end{equation*}
This completes the proof.
\qed

\begin{lem}\label{key lemma: Witt-Burnside2}
Let $G$ be an abelian profinite group.
Then, for every open subgroup $U$ of $G$, it holds that
$$Y_{U}-X_{U}\in
\mathbb Z\,\left[\,\frac {1}{c_i},
\frac {1}{d_j} \,:\, 1\le i \le s, 1\le j \le t\,\right]
\left[X_{S}\,:\,U\subsetneqq S \subseteq G, \,\,\,{\rm and }\,\, S\in \mathcal O(G)\,\right].$$
\end{lem}
\noindent{\bf Proof.}
Let us first assume that $q$ and $r$ are nonzero.
Note that
$Y_{G}=X_{G}$.
Now, we assume that
$$Y_{V}-X_{V}\in
\left(\mathbb Z_G \cap \mathbb Z\left[\frac 1q,\frac 1r\right]\right)
\left[X_{S}\,:\, V\subsetneqq S \subseteq G, \,\,\,{\rm and }\,\, S\in \mathcal O(G)\,\right]$$
for all $V \in \mathcal O(G)$ with $U \subsetneqq V \subseteq G.$
Since $G$ is abelian, Eq. \eqref{r-phi=q-phi:Witt-Burnside}
is transformed into
\begin{equation*}\label{r-phi=q-phi:Witt-Burnside:transform}
\begin{aligned}
Y_{U}-X_{U}=\sum_{V \in \mathcal O(G) \atop U \subsetneqq V \subseteq G}
\frac {q^{(V:U)-1} X_{V}^{(V:U)}-r^{(V:U)-1}Y_{V}^{(V:U)}}
{(V:U)}.
\end{aligned}
\end{equation*}
Note that if a prime $p$ divides $(V:U)$ and $q$  then
it cannot be a divisor of the denominator of the irreducible fraction of
${q^{(V:U)-1}}/{(V:U)}.$
Also, it cannot be a divisor of the denominator of the irreducible fraction of
${r^{(V:U)-1}}/{(V:U)}$ if it divides $r$.
Now, our assertion follows from the induction hypothesis.

Next, we assume that $q$ is zero and $r$ is nonzero.
In this case,
Eq. \eqref{r-phi=q-phi:Witt-Burnside} is reduced to
\begin{equation*}
\begin{aligned}
Y_{U}-X_{U}=-\sum_{V \in \mathcal O(G) \atop U \subsetneqq V \subseteq G}
\frac {r^{(V:U)-1}}
{(V:U)}Y_{V}^{(V:U)}.
\end{aligned}
\end{equation*}
Since a prime $p$ dividing $(V:U)$ and $r$
cannot be a divisor of the denominator of the irreducible fraction of
${r^{(V:U)-1}}/{(V:U)}$ we obtain the desired result by induction hypothesis.
So we are done.
\qed
\vskip 3mm
{\bf Proof of Theorem \ref{main theorem:Witt- Burnside}.}
The ``if" part follows from Lemma \ref{key lemma: Witt-Burnside2}.
For the ``only if" part,
let us assume that there exists a unique strict-isomorphism,
say $\tau_q^r: \mathbb W_G^q(A) \to \mathbb W_G^r(A)$.
Hence, we have
$\Phi_G^q=\Phi_G^r\circ \tau_q^r$.
Assume that a prime $p$ divides $q$ and $(G:U)$
for some open subgroup $U$ of $G$,
but not $r$.
Consider a maximal filtration
$G=V_1\supsetneqq V_2 \supsetneqq \cdots \supsetneqq V_k=U$,
where $V_i\, (1\le i \le k)$ are open subgroups of $G$
such that there is no open subgroup between $V_i$ and $V_{i+1}$.
Then $p$ must divide $(V_i:V_{i+1})$ for some $1\le i \le k$.
Letting
\begin{equation}\label{unity condi: Witt-Burnside}
(a_{U})_{U\in \mathcal O(G)}
=\tau_q^r(0,0,\cdots, \stackrel {V_i-{\rm th}}{\overbrace {1}},0,0,\cdots),
\end{equation}
we can deduce the identity
\begin{equation}\label{unity condi: Witt-Burnside1}
q^{(V_{i}:V_{i+1})-1}-r^{(V_{i}:V_{i+1})-1}
=(V_i:V_{i+1})a_{V_{i+1}}
\end{equation}
by comparing the $V_{i+1}$-th component of both sides
of Eq. \eqref{unity condi: Witt-Burnside}.
We claim that $p$ is a unit in $A$.
To see this, transform Eq. \eqref{unity condi: Witt-Burnside1} into
\begin{equation*}\label{unity condi: Witt-Burnside2}
-r^{(V_{i}:V_{i+1})-1}
=p\left(\frac {(V_i:V_{i+1})a_{V_{i+1}}-q^{(V_{i}:V_{i+1})-1}}{p}\right)
\end{equation*}
Since $p$ and $r^{(V_{i}:V_{i+1})-1}$ are coprime
there are $x,y \in \mathbb Z$ such that
$$px-r^{(V_{i}:V_{i+1})-1}y=1.$$
Therefore,
$$-r^{(V_{i}:V_{i+1})-1}y=1-px=
p\left(\frac {(V_i:V_{i+1})a_{V_{i+1}}-q^{(V_{i}:V_{i+1})-1}}{p}\right)y.$$
Hence,
$p\left(x+
\left(\frac {(V_i:V_{i+1})a_{[V_{i+1}]}-q^{(V_{i}:V_{i+1})-1}}{p}\right)y\right)=1.$
This justifies our claim.
Similarly, a prime $p$ dividing $r$ and $(G:U)$, but not $q$,
should be a unit.
This implies that $A$ should be a
$\mathbb Z\,[\,\frac {1}{c_i},
\frac {1}{d_j} \,:\, 1\le i \le s, 1\le j \le t\,]$-algebra.
\qed
\vskip 3mm
From Theorem \ref{main theorem:Witt- Burnside} it is immediate that
$\mathbb W_G^q(\mathbb Z)$ is classified up to strict-isomorphism
by the set of prime divisors of $q$ contained in $D^{\rm pr}(G)$.
Thus, we proved
\begin{cor}
Let $q$ vary over the set of integers and $G$ be an abelian profinite group.
Then, $\mathbb W_G^q$ is classified up to strict natural isomorphism
by the set of prime divisors of $q$ contained in $D^{\rm pr}(G)$.
\end{cor}

\vskip 3mm

{\bf Acknowledgement} \ \ \
The author would like to express his sincere gratitude
to the referee for his/her correction of many errors
in the previous version and valuable advices.
\small{

}
\end{document}